\documentclass{ws-ijbc}
\usepackage{graphicx}
\usepackage{epstopdf}
\usepackage{float}
\usepackage{subfigure}
\usepackage{color}
\usepackage{hyperref}
\usepackage{multicol}

\begin{document}

\catchline{}{}{}{}{} % Publisher's Area please ignore

\markboth{Danielle F. P. Toupo, David G. Rand, Steven H. Strogatz}{REPLICATOR-MUTATOR MODEL OF PRISONER'S DILEMMA}

\title{LIMIT CYCLES SPARKED BY MUTATION \\ IN THE  REPEATED PRISONER'S DILEMMA}
%USING \LaTeX\footnote{For the title, try not to use more than
%three lines. Typeset the title in 15 pt Times Roman, uppercase and
%boldface.}}
\author{DANIELLE F. P. TOUPO}
%\footnote{Typeset names in 11 pt Times Roman.
%Use the footnote to indicate the present or permanent address of
%the author.}}

\address{Center for Applied Mathematics, Cornell University\\
Ithaca, NY 14853, USA\\
dpt35@cornell.edu}
%\footnote{State completely without
%abbreviations the affiliation and mailing address, including
%country. Typeset in 11~pt Times Italic.}}

\author{DAVID G. RAND}
\address{Department of Psychology, Yale University\\
New Haven, CT 06511, USA\\
david.rand@yale.edu}

\author{STEVEN H. STROGATZ}
\address{Department of Mathematics, Cornell University\\
Ithaca, NY 14853, USA\\
strogatz@cornell.edu}

\maketitle

\begin{history}
\received{(to be inserted by publisher)}
\end{history}

\begin{abstract}

We explore a replicator-mutator model of the repeated Prisoner's Dilemma involving three strategies: always cooperate (ALLC), always defect (ALLD), and tit-for-tat (TFT). The dynamics resulting from single unidirectional mutations are considered, with detailed results presented for the mutations TFT $\rightarrow$ ALLC and ALLD $\rightarrow$ ALLC. For certain combinations of parameters, given by the mutation rate $\mu$ and the complexity cost $c$ of playing tit-for-tat, we find that the population settles into limit cycle oscillations, with the relative abundance of ALLC, ALLD, and TFT cycling periodically.  Surprisingly, these oscillations can occur for unidirectional mutations between \emph{any} two strategies. In each case, the limit cycles are created and destroyed by supercritical Hopf and homoclinic bifurcations, organized by a Bogdanov-Takens bifurcation. Our results suggest that stable oscillations are a robust aspect of a world of ALLC, ALLD, and costly TFT; the existence of cycles does not depend on the details of assumptions of how mutation is implemented.
\end{abstract}
{%begin{multicols}{2}
\section{Introduction}
\noindent 
Cooperation, where individuals pay costs to benefit others, is a cornerstone of human civilization. By cooperating, people create value and thus increase the ``size of the pie.'' This makes a group in which everyone cooperates better off than a group where everyone is selfish. Cooperation can be difficult to achieve, however, because creating that collective benefit is often individually costly. How, then, could the selfish process of natural selection give rise to such altruistic cooperation? Evolutionary game theorists have devoted a great deal of effort to answering this question \cite{Axelrod84, Boyd03, Boyd92, Fowler05, Fu08, Fudenberg90,  Hamilton64, Hauert02, Helbing09, Manapat12, May87, McNamara08, Nakamaru97, Nowak06, Nowak04, Panchanathan04, Rand12, Rand131, Rand132, Szolnoki09, Szolnoki04, Tarnita09, Traulsen06, Trivers71, Wedekind00}. 

The standard game-theoretic paradigm for studying cooperation is the Prisoner's Dilemma \cite{Rapoport65}. Two players simultaneously choose to either cooperate (C) or defect (D), and each receives a payoff depending on the two choices. Two cooperators both earn the reward of mutual aid $R$, while two defectors receive a punishment $P$. If one player cooperates while the other defects, the defector earns the temptation payoff $T$ while the cooperator receives the sucker payoff $S$. If the relationship $T > R > P > S$ holds, the game is a Prisoner's Dilemma because mutual cooperation is better than mutual defection ($R > P$), but individually defectors always out-earn cooperators ($T>R$ if the partner cooperates, and $P > S$ if the partner defects). Thus the Prisoner's Dilemma captures the tension between individual and collective interests, the conundrum at the heart of cooperation. 

To study the evolution of cooperation, evolutionary game theorists typically combine game theory with differential equations to create an evolutionary dynamic. The replicator equation is one of the most common such models \cite{Hofbauer79}: strategies with above-average payoffs become more common over time while strategies with below-average payoffs become less common. As described above, defectors always out-earn cooperators. Thus in a simple world of pure cooperators versus pure defectors, evolution via the replicator equation always leads to the extinction of cooperation and a population made up solely of defectors.

How, then, do we explain the success of cooperation which is so evident in the world around us? Numerous mechanisms for the evolution of cooperation have been proposed \cite{Nowak06}, and empirical evidence has been provided for their importance in human cooperation \cite{Rand131}. Chief among these mechanisms is direct reciprocity, also known as reciprocal altruism \cite{Axelrod84, Nowak92, Nowak93, Rand09, Trivers71}: when agents interact repeatedly, evolution can favor cooperation. If I will only cooperate with you in the next period if you cooperate with me in the current period, cooperation can be the payoff-maximizing strategy (as long as the game continues to the next period with high enough probability). 

Tit-for-tat (TFT) is the most well-known of these reciprocal strategies. TFT begins by cooperating, and then merely copies its opponent's move in the previous period. Thus cooperators receive cooperation and profit, whereas defectors receive defection and lose out. Moreover, a population in which everyone plays TFT is resistant to invasion by defectors. If a lone defector is introduced into  such a population, that ALLD player receives a low payoff relative to the resident TFT players. As a result, selection disfavors the invader and TFT is evolutionarily stable: repeated interactions promote the evolution of cooperation. 

Tit-for-tat, however, has an Achilles' heel \cite{Nowak06}: in some situations it can be invaded by kinder, gentler strategies \cite{Boyd87, Imhof05, van10, van12}. To see this, imagine again a population where everyone plays TFT, except for one player who always cooperates (ALLC). Both the resident TFT players and the ALLC deviant cooperate in every round, and thus receive equal payoffs. This means that neutral drift (in the sense of population genetics) can allow ALLC to increase in frequency, a process called neutral invasion. A more serious weakness of the TFT players is that other factors can impose costs on them -- costs which ALLC players avoid. For example, TFT is a more sophisticated strategy than ALLC, and hence may incur complexity costs. These costs arise because TFT needs to spend energy interpreting the other player's last move before it can respond appropriately, whereas ALLC is an unconditional strategy. On top of that, the vindictiveness of TFT can hurt it in noisy or error-prone environments. If players sometimes make mistakes and defect when they meant to cooperate,  such a mistake would go unnoticed by ALLC, but would send two TFT players into a vendetta of retaliatory defections. 

In the presence of such costs, the TFT residents are ultimately overtaken by the ALLC invaders. And there's the rub:  once ALLC becomes sufficiently common at TFT's expense, it opens the door to invasion by nasty, uncooperative ALLD players, who can ruthlessly exploit ALLC and sweep through the population. Through this sort of scenario, TFT populations can eventually wind up succumbing to defection, suggesting that cooperation may be doomed even in repeated games. 

Here we present a solution to this problem: we show that incorporating mutation into the replicator equation breaks ALLD's dominance over the evolutionary outcomes by giving rise to stable cycles involving substantial levels of cooperation. Surprisingly, it is not necessary for every strategy to mutate into every other strategy in order to get cyclical behavior, or even for ALLD to mutate into TFT. We find that adding \emph{any}  single unidirectional pathway for mutation can lead to stable cycles. Thus we provide evidence that stable oscillations are a robust aspect of a world of ALLD, ALLC and costly TFT players, and that direct reciprocity can lead to substantial cooperation even in the face of invasion by unconditional cooperators.

\subsection{Relation to previous work}
To allow for mutation, the relevant mathematical setting needs to change from the replicator equation to the replicator-mutator equation \cite{Nowak06, Stadler92}. A number of previous authors have studied limit cycles in replicator-mutator equations, motivated by applications to language change \cite{Mitchener04}, autocatalytic chemical reaction networks \cite{Stadler92}, evolutionary games \cite{Bladon10, Galla11, Imhof05}, and multi-agent decision making \cite{Pais11, Pais12}. 

Our work is most closely related to that of Imhof \emph{et al.} [2005]. They studied the evolutionary game dynamics of ALLD, ALLC, and TFT for finite populations, where stochastic effects become important. One of their most striking results was the observation of the phenomenon mentioned above -- an evolutionary cycle that goes from ALLD to TFT to ALLC and back to ALLD again. Even more remarkably, they found that the cycle spends nearly all its time lingering in the vicinity of TFT, even though ALLD is a strict Nash equilibrium! 

In the model considered by Imhof \emph{et al.} [2005], the conditional strategy TFT was assumed to incur a complexity cost $c$, relative to the simpler unconditional strategies ALLC and ALLD. Mutation was assumed to be uniform: every strategy mutates to the other two with equal probability. More precisely, if $i \neq j$, strategy $i$ mutates to strategy $j$ with probability $\mu$ and stays the same with complementary probability $1-2\mu$. 

For the deterministic case of infinitely large populations, Imhof \emph{et al.} [2005] found that for certain combinations of parameters $c$ and $\mu$, the replicator-mutator equations have a stable limit cycle, corresponding to the evolutionary cycles observed in their simulations. They mentioned bifurcations associated with the birth and death of the limit cycle, but did not present a bifurcation analysis or a stability diagram to locate the bifurcation curves in the $(\mu, c)$ parameter space. 

We were curious to learn more about the evolutionary cycle seen in the infinite population model. What bifurcations create and destroy this limit cycle? How does its bifurcation structure depend on the details of how the mutations are implemented, in a graph-theoretic sense? That is, if we think of the three strategies as the vertices of a triangle graph, with mutations occurring along the edges between them, the uniform mutation case studied by Imhof \emph{et al.} [2005] amounts to a complete graph with equal weights $\mu$ on its six directed edges. What would happen, by contrast, in the opposite extreme case of unidirectional mutation along one of these six directed edges? 

We found that for all six possible unidirectional mutations, a stable limit cycle exists in a certain part of $\left(\mu, c\right)$ space. The cycle always oscillates in the same rotational sense, moving from the neighborhood of ALLD to TFT to ALLC and back toward ALLD again, just as it does in the case of uniform mutation. Moreover, the commonalities extend to the types and locations of the bifurcations that create and destroy the cycle. The region in parameter space where stable limit cycles exist is always bounded on one side by a curve of supercritical Hopf bifurcations and on the other side by a curve of homoclinic bifurcations. These two curves meet tangentially at a Bogdanov-Takens point at one end of the stability region for the limit cycle. All of these statements are true of the uniform mutation case as well. 

This paper is organized as follows. In the next section we review the formulation of the Prisoner's Dilemma and its associated replicator equation, followed by their generalization to the replicator-mutator equation. Then we focus on two special cases of unidirectional mutation, TFT~$\rightarrow$~ALLC and ALLD~$\rightarrow$~ALLC, and summarize the results from the other four cases, as well as for the case of uniform mutation. The paper concludes with a conjecture and a brief discussion.

\section{Model}

\subsection{Prisoner's Dilemma}
Following Axelrod [1984], we fix the parameter values $\displaystyle T = 5, R = 3,  P = 1, S = 0$. These satisfy the inequalities  required for the game to qualify as a Prisoner's Dilemma: $\displaystyle T > R > P > S$ and $\ R > \left(T + S\right)/2$. The final inequality implies that if the two players play many rounds with one another, it is better for both of them to cooperate all the time rather than engage in alternating bouts of getting suckered and suckering the other.  

Now consider a repeated Prisoner's Dilemma among players using the strategies ALLC, ALLD, and TFT.  Then, in the limit where the players meet infinitely often and there is no discounting of future interactions, their average payoffs are given by the  payoff matrix given in Table 1 (the entries show the average payoff that the row player gets when playing against the specified column player):

\begin{table}[h]
\tbl{Payoff matrix of repeated prisoner's dilemma}
{\begin{tabular}{l c c c c c}\\[-2pt]
\toprule
{} &{} &ALLD &TFT &ALLC\\[6pt]
\hline\\[-2pt]
{}&ALLD &$P$ &$P$ &\phantom0$T$\\[1pt]
{}&TFT &$P$ &$R$ &\phantom0$R$\\[2pt]
{}&ALLC &$S$ &$R$ &\phantom0$R$\\[1pt]
\botrule
\end{tabular}}
\end{table}

\noindent For example, an ALLC player gets suckered every time against an ALLD player, and therefore receives an average payoff of $S$ in their depressing encounters. But when TFT plays ALLD, it gets suckered only on the first round, but after that it reciprocates each defection with defection, leading to an infinitely long string of mutual defections and hence an average payoff of $P$.   

\subsection{Replicator equation}
Next we set the game in an evolutionary framework. If each strategy reproduces at a rate proportional to its relative fitness, the resulting dynamics can be approximated by the following set of ordinary differential equations, known as the replicator equations:
\begin{eqnarray}
\dot{x} &=& x \left(f_x - \phi \right) \nonumber\\
\dot{y} &=& y \left(f_y - \phi \right) \nonumber\\
\dot{z} &=& z \left(f_z - \phi \right).  \label{eqn:replicator}
\end{eqnarray}
Here $x, y$ and $z$ denote the fractions of the population playing ALLD, TFT, and ALLC, respectively; 
$f_i$ is the fitness of strategy $i$, defined as its expected payoff against the current mix of strategies; 
and 
\begin{eqnarray}
\phi = x f_x + y f_y + z f_z 
\end{eqnarray}
is the average fitness in the whole population. 
By summing the differential equations for $\dot{x}, \dot{y},$ and $\dot{z}$, one can verify that 
$\displaystyle x + y + z = 1$ for all time, as required by the definition of $x, y$ and $z$ as relative frequencies. 

The payoff matrix for the repeated Prisoner's Dilemma implies that the fitnesses are given by 
\begin{eqnarray}
f_x &=& x P+y P+z  T \nonumber\\
f_y &=& x P+y  R+z  R \nonumber\\
f_z &=& x S+y  R+z  R. \label{eqn:fitness}
\end{eqnarray}
For the parameter values $T = 5, R = 3,  P = 1$, and $S = 0$ assumed above, and by replacing $z$ with $1-x-y$, 
we find that the fitnesses reduce to
\begin{eqnarray}
f_x &=& 5 -4  x -4  y   \nonumber\\
f_y &=& 3 - 2  x \nonumber\\
f_z  &=& 3 - 3  x \label{eqn:fitness_with_numbers}
\end{eqnarray}
and the average fitness in the population becomes $\displaystyle \phi = 3 - x\left(1 +x+3  y\right).$

The phase portrait for the replicator equation \eqref{eqn:replicator} with fitnesses \eqref{eqn:fitness_with_numbers} can be drawn in the $(x,y)$ plane, by eliminating $z$ via $z=1-x-y$. Unfortunately this way of presenting the phase portrait has certain disadvantages. It distorts the geometry of the trajectories and it arbitrarily highlights two of the strategies at the expense of the third. 
In the original $(x,y,z)$ phase space, the phase portrait lives on the equilateral triangle defined by the face of the simplex $x+y+z=1$, where $0 \leq x, y, z \leq 1$. Thus, a more appealing and symmetrical approach is to show the phase portrait as it actually appears on simplex. The following change of variables achieves this goal:
\begin{equation}
 \binom{X}{Y} =  \begin{pmatrix} 1 & \frac{1}{2}  \\ 0 & \frac{\sqrt{3}}{2} \end{pmatrix} \binom{x}{y} \label{eqn:transformation}
\end{equation}
Equation \eqref{eqn:transformation} can be shown to be equivalent, up to a uniform scaling and a translation, to the change of variables given by Eq. 35 in Wesson \& Rand [2013].

Figure~\ref{fig:c0mu0} plots the resulting phase portrait. The first thing to note is that it contains a saddle point at $(x,y,z)=(1,0,0)$, corresponding to the entire population at ALLD. The fact that the ALLD corner is a saddle point, rather than a stable fixed point, reflects our simplifying assumptions that the game is infinitely repeated with no discounting of future interactions; in effect, the first round of the game is totally ignored. Without these assumptions, a standard result in repeated games \cite{Nowak06} is that there is bistability along the ALLD-TFT edge. In that case the ALLD corner would have a basin of attraction whose size depends on the payoffs. 

\begin{figure}[!htp]
\centering
\includegraphics[scale=0.2]{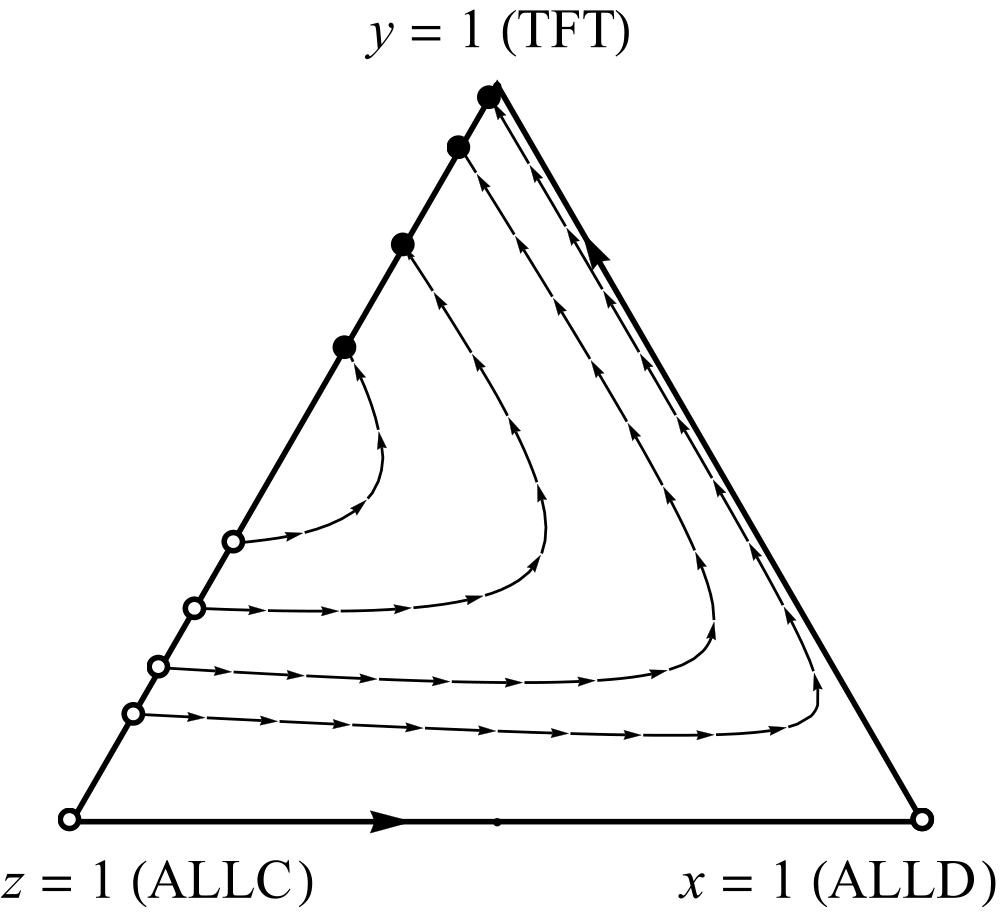} %case mu = 0 c= 0
\caption{Phase portrait of system \eqref{eqn:replicator} with $c=0$ and $\mu=0$.}
\label{fig:c0mu0}
\end{figure}

Second, observe that Fig.~\ref{fig:c0mu0} displays a line of neutrally stable fixed points along the side joining ALLC to TFT, where ALLD is absent. Neutral drift takes place along this side in finite populations. In the infinite population case shown here, the system almost always ends up in a nirvana of cooperation, with a mix of TFT and ALLC determined by the initial conditions, and with ALLD extinct. 

However, the line of neutrally stable fixed points in this phase portrait is a structurally unstable feature. If the governing equations are perturbed by the addition of arbitrarily small terms, one expects that the line of fixed points will break and be replaced by something qualitatively different.

Indeed, when we associate a small complexity cost $c$ to playing TFT, the payoff matrix changes to that shown in Table 2:

\begin{table}[!h]
\tbl{Payoff matrix of replicator model with cost.}
{\begin{tabular}{l c c c c c}\\[-2pt]
\toprule
{} &{} &ALLD &TFT &ALLC\\[6pt]
\hline\\[-2pt]
{}&ALLD &$P$ &$P$ &\phantom0$T$\\[1pt]
{}&TFT &$P - c$ &$R - c$ &\phantom0$R - c$\\[2pt]
{}&ALLC &$S$ &$R$ &\phantom0$R$\\[1pt]
\botrule
\end{tabular}}
\end{table}

The new fitnesses become 
\begin{eqnarray}
f_x &=& 5 - 4  x - 4  y   \nonumber\\
f_y &=& 3 - c - 2  x \nonumber\\
f_z  &=& 3 - 3  x \label{eqn:fitness_with_numbers_and_cost}
\end{eqnarray}
and 
\begin{eqnarray}
\phi = 3 - c  y - x (1 + x + 3  y). \label{eqn:avg_fitness_with_numbers_and_cost}
\end{eqnarray}
In the corresponding phase portraits, shown in Fig.~\ref{fig:c>0mu0}, the previous line of neutrally stable fixed points turns into an invariant line with the vector field flowing from a saddle point at TFT to another saddle point at ALLC.\\
\begin{figure}[h]
\centering
\subfigure[\  $\displaystyle 0 < c < \frac{2}{3} $]{%
\includegraphics[scale = 0.15]{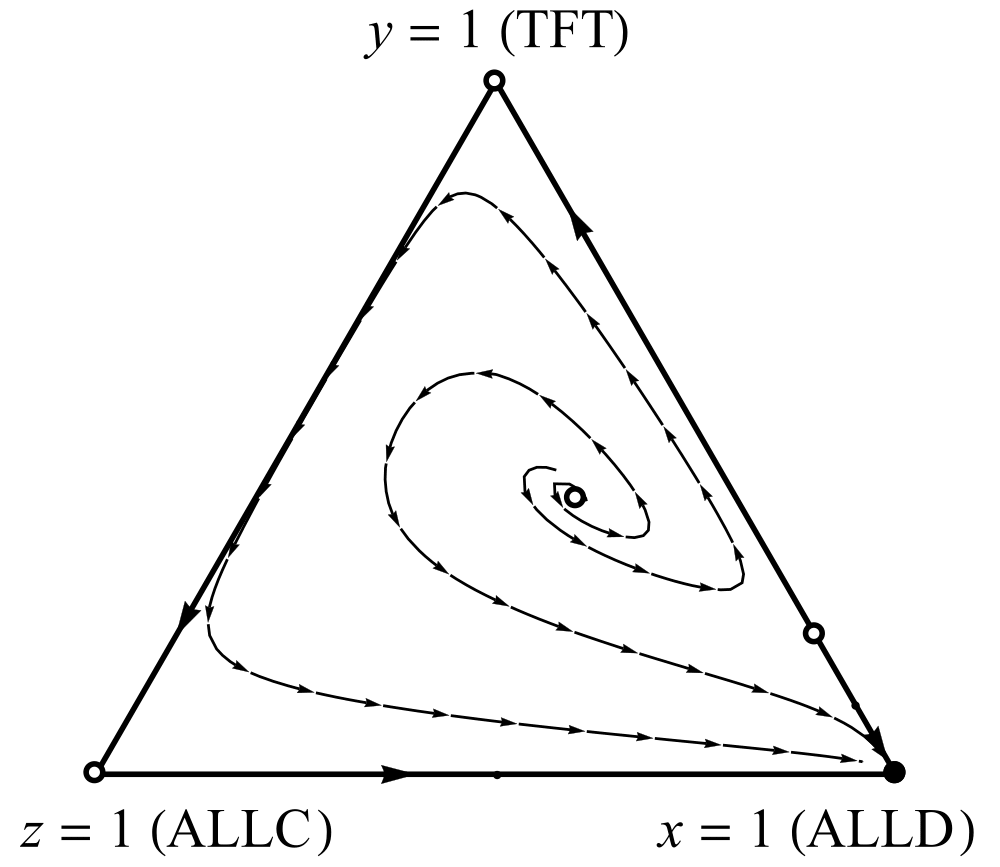}}
%\label{fig:subfigure1}}
\quad
\subfigure[\  $ \displaystyle \frac{2}{3} \leq c < 2 $]{%
\includegraphics[scale = 0.15]{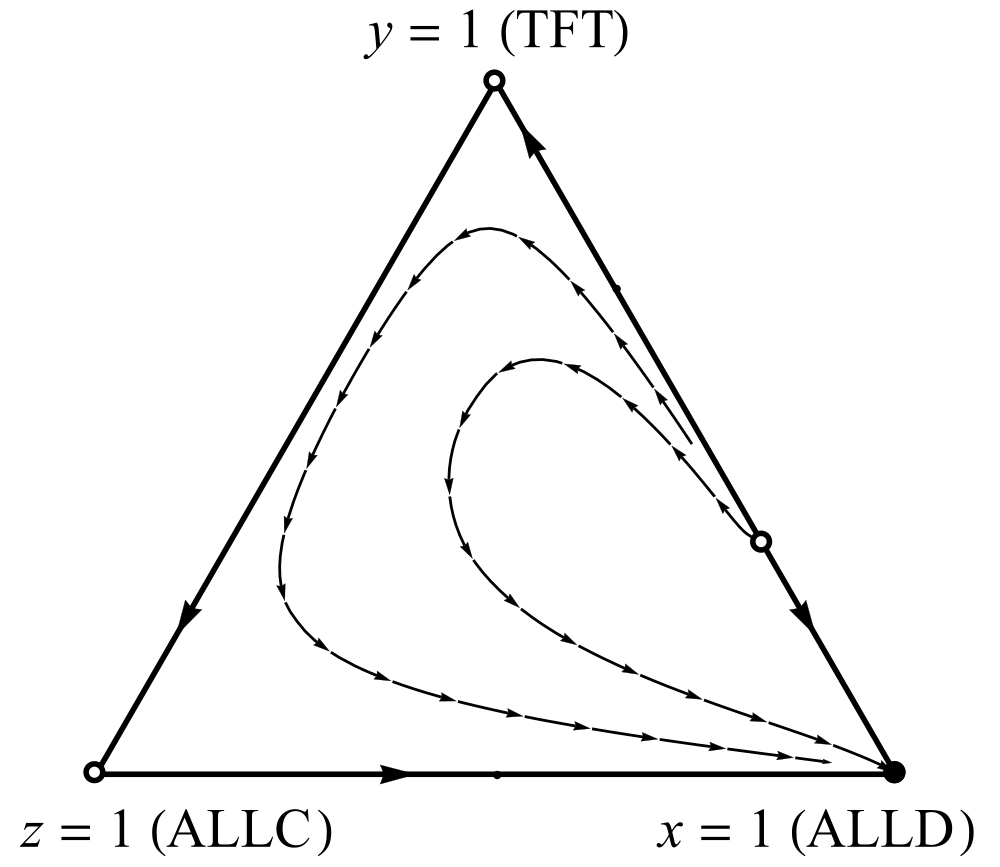}}
%\label{fig:subfigure2}}
\quad
\subfigure[\  $\displaystyle c \ge 2$]{%
\includegraphics[scale = 0.15]{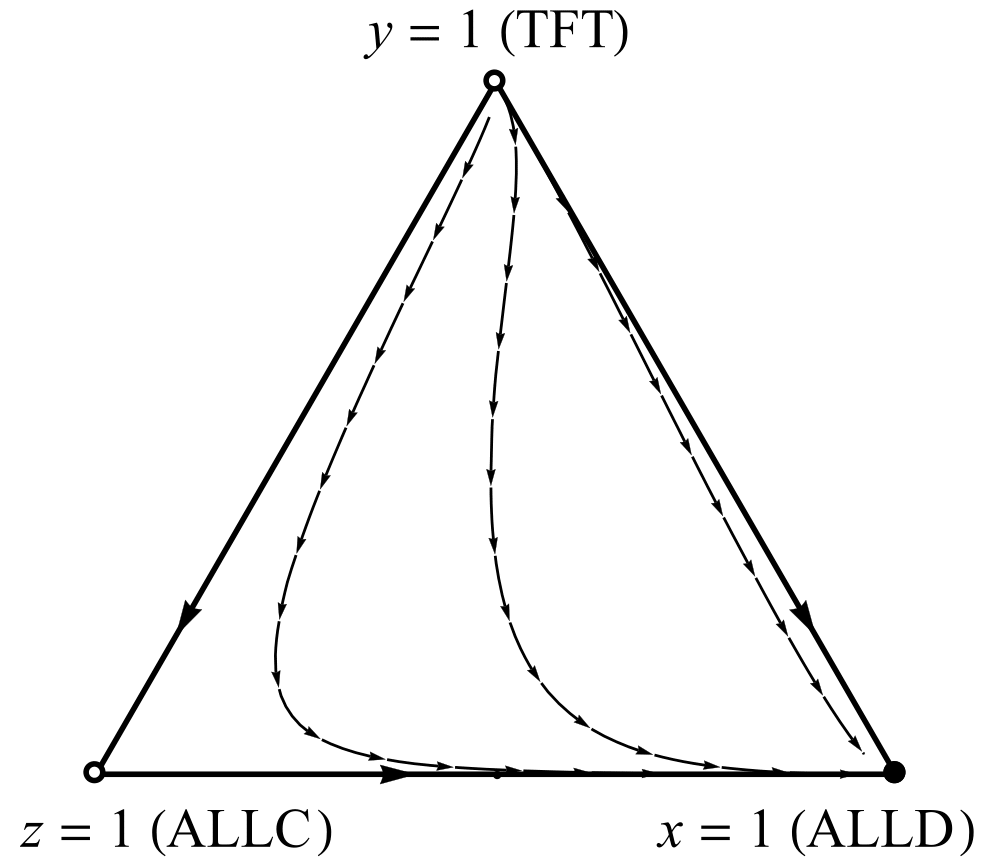}}

\caption{Phase portrait of system \eqref{eqn:replicator} with $c > 0$ and  $\mu=0$.}
\label{fig:c>0mu0}
\end{figure}

\noindent The structure of the rest of the phase portrait depends on the size of $c$, as shown in Fig.~\ref{fig:c>0mu0}, but the long-term behavior of the system is the same in all three cases: ALLD takes over and cooperation dies out.

Notice that none of the phase portraits so far contain any limit cycles.  This could have been anticipated from a general theorem that  forbids limit cycles in any system of replicator equations involving $n=3$ strategies, for any game and any payoff matrix \cite{Hofbauer98}.  Periodic solutions can exist, but only within continuous families of neutrally stable cycles. Such periodic orbits are not isolated and hence do not qualify as limit cycles.

 \section{Replicator-Mutator Equation}
 Limit cycles do, however, become possible when we allow ALLD, ALLC, and TFT to mutate into one another. For simplicity, let us restrict attention to single  unidirectional mutations only. Then there are six possibilities.
 
 %------------------------------  TFT -> ALLC------------------------------------------------------------------------------
\subsection{Example 1: TFT $\rightarrow$ ALLC}
In this case, we assume that after replication occurs, a player with strategy TFT mutates to ALLC with probability $\mu$. 
Then the replicator-mutator system is 
\begin{eqnarray}
\dot{x} &=& x\left(f_x - \phi\right)  \nonumber\\
\dot{y} &=& y   \left(f_y - \phi\right) - \mu  y  f_y.   \label{eqn:rmeqn1_TFTALLC}
\end{eqnarray}
where $x, y,$ and $z$ denote the frequencies of ALLD, TFT, and ALLC, respectively. Note that as before, we have chosen to eliminate $z$ from the equations by using the identity $z=1-x-y$. Likewise, we do not write the equation for $\dot{z}$ explicitly, since it can be obtained from $\dot{x}$ and $\dot{y}$ if needed via $\dot{z} = -\dot{x}-\dot{y}$.

After insertion of  \eqref{eqn:fitness_with_numbers_and_cost} and  \eqref{eqn:avg_fitness_with_numbers_and_cost} into \eqref{eqn:rmeqn1_TFTALLC}, the replicator-mutator system becomes
\begin{eqnarray}
\dot{x} &=& x  \left[ \left(c-4\right)   y+x  \left(x+3   y-3\right)+2 \right], \nonumber\\
\dot{y} &=& y   \left[ c   \left(\mu +y-1\right)-3   \mu +x \left(2   \mu +x+3   y-1\right) \right].  \label{eqn:rmeqn2_TFTALLC}
\end{eqnarray}
The right hand side is cubic in $x$ and $y$. With the help of computer algebra, one can calculate the fixed points for the system, all but one of which lie on the boundary of the simplex (the equilateral triangle). Explicit formulas for these fixed points are presented in  Appendix A. We have also calculated the curves in $(\mu, c)$ space at which the interior fixed point undergoes saddle-node and supercritical Hopf bifurcations; see Appendix A for details. As the parameters are varied, the stable limit cycles that are created by the supercritical Hopf bifurcation are later destroyed by a homoclinic bifurcation. The associated curve of homoclinic bifurcations was computed numerically with the help of the continuation package MATCONT \cite{Govaerts08}. 

Figure~\ref{fig:bifnTFTALLC} plots the bifurcation curves in $(\mu, c)$ space. The saddle-node curve is shown in blue, the Hopf curve in red, and the homoclinic curve in green. These curves partition the parameter space into four regions, marked 1, 2, 3, and 4 on the figure. 

\begin{figure}[!h]
\centering
\includegraphics[scale = 0.27]{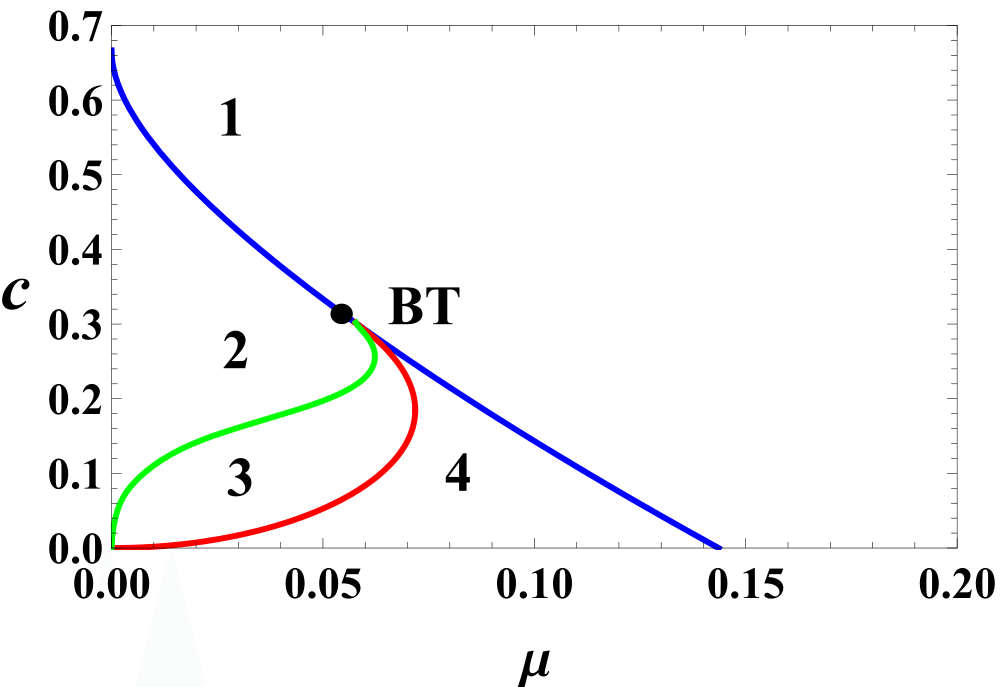}
\caption{Stability diagram in $(\mu, c)$ space when TFT $\rightarrow$ ALLC. Saddle-node curve, blue; supercritical Hopf curve, red; homoclinic curve, green; BT, Bogdanov-Takens point. } 
\label{fig:bifnTFTALLC}
\end{figure}

Figure~\ref{fig:TFTALLCphaseportraits} shows representative phase portraits for each of the four regions. In region~1, ALLD wipes out the other two strategies. (This makes sense intuitively. In region 1, TFT pays a high cost in fitness compared to the other two strategies. Furthermore, because $\mu$ is large in region 1, TFT mutates rapidly into ALLC, which in turn is clobbered by ALLD.) In region~2,  a new pair of fixed points exist; they were created in a saddle-node bifurcation on the right-hand boundary of the equilateral triangle when the parameters crossed the saddle-node bifurcation curve. The unstable spiral seen in region~2 is the descendant of that node. Despite the existence of these two new fixed points, the system's long-term behavior remains the same as in region~1: almost all trajectories approach ALLD.  

\begin{figure}[h]
\centering
\subfigure[\  Region 1]{%
\includegraphics[scale = 0.13]{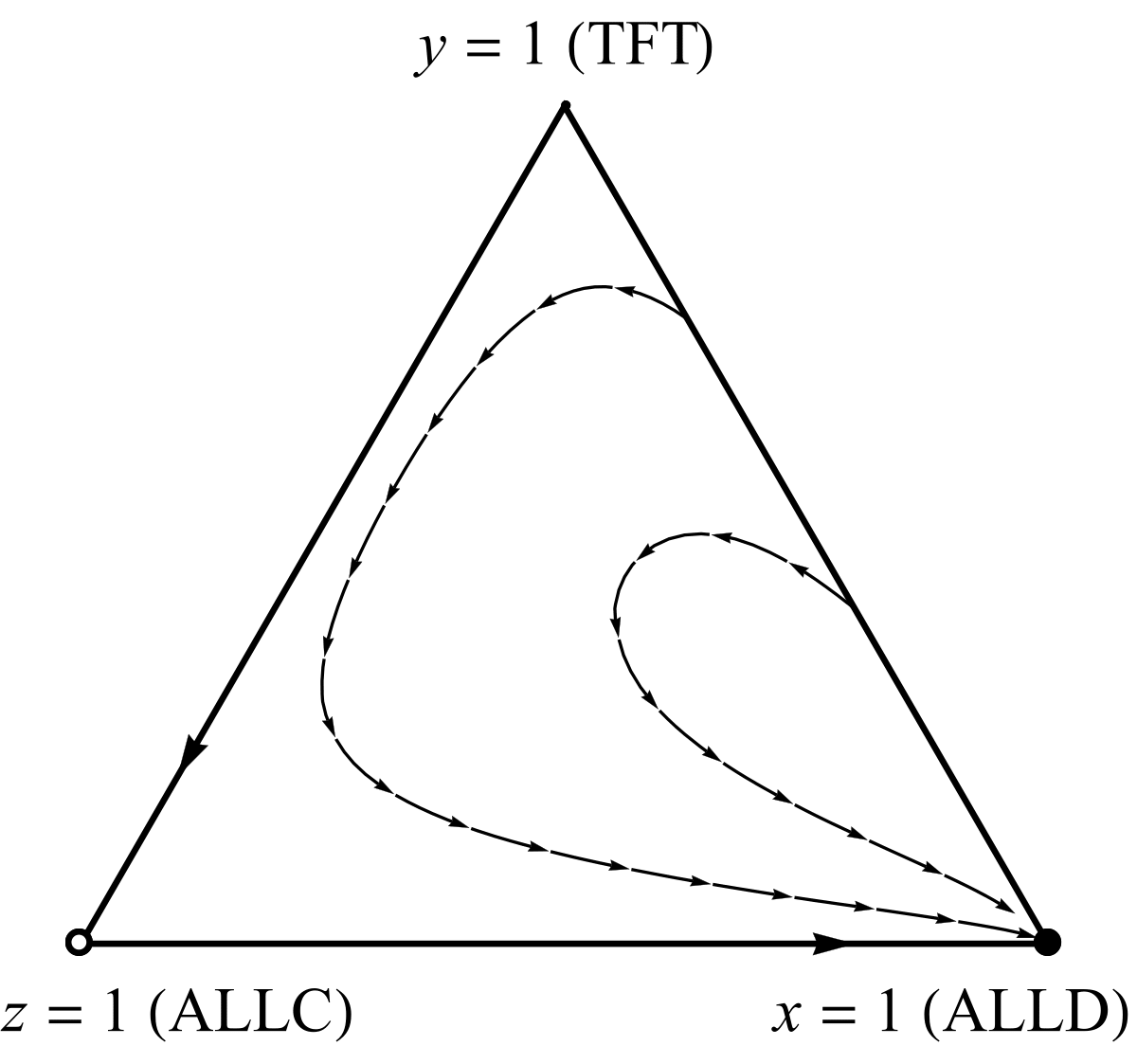}}
%\label{fig:subfigure1}}
\quad
\subfigure[\  Region 2]{%
\includegraphics[scale = 0.13]{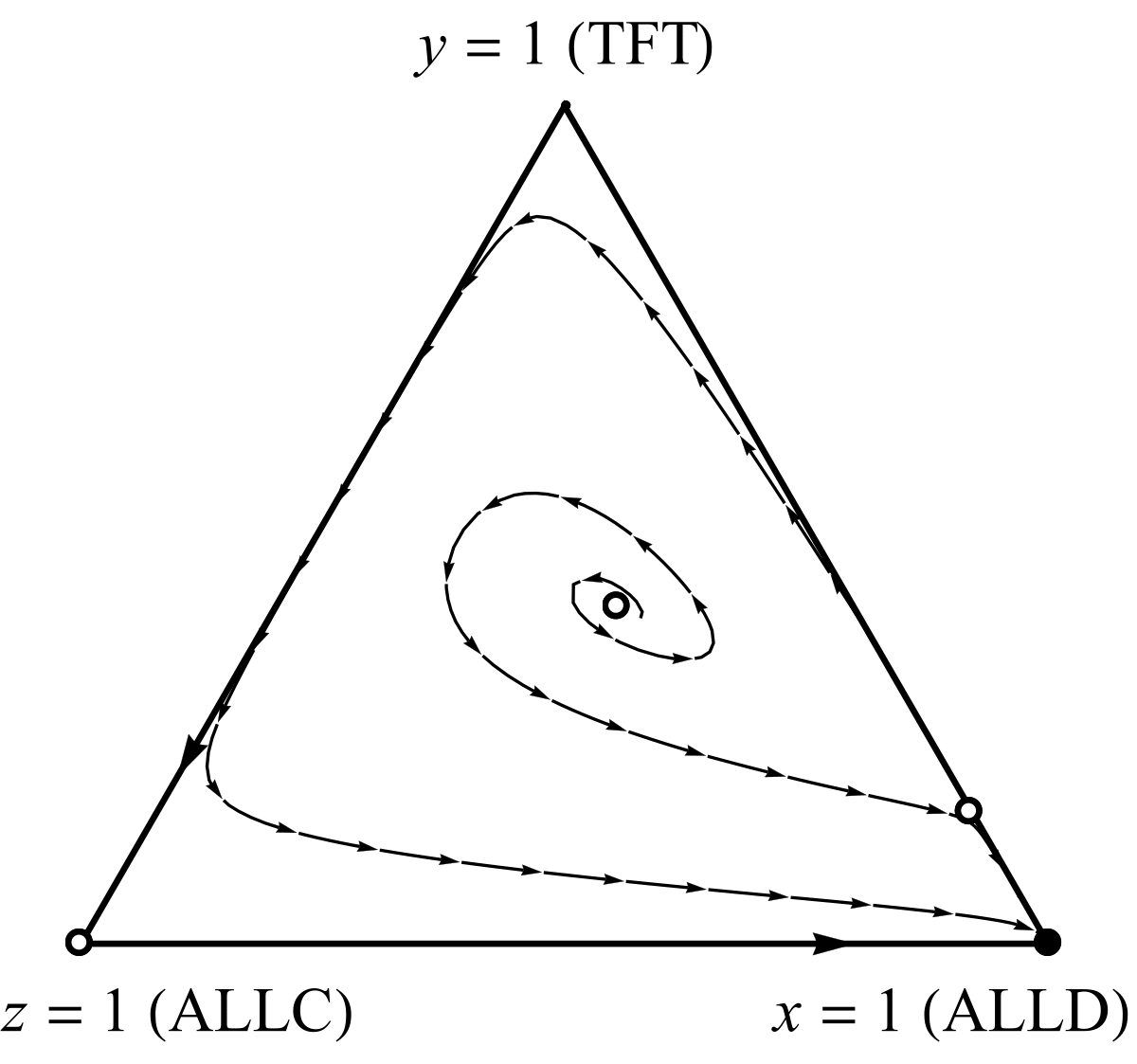}}
%\label{fig:subfigure2}}

\subfigure[\  Region 3]{%
\includegraphics[scale = 0.13]{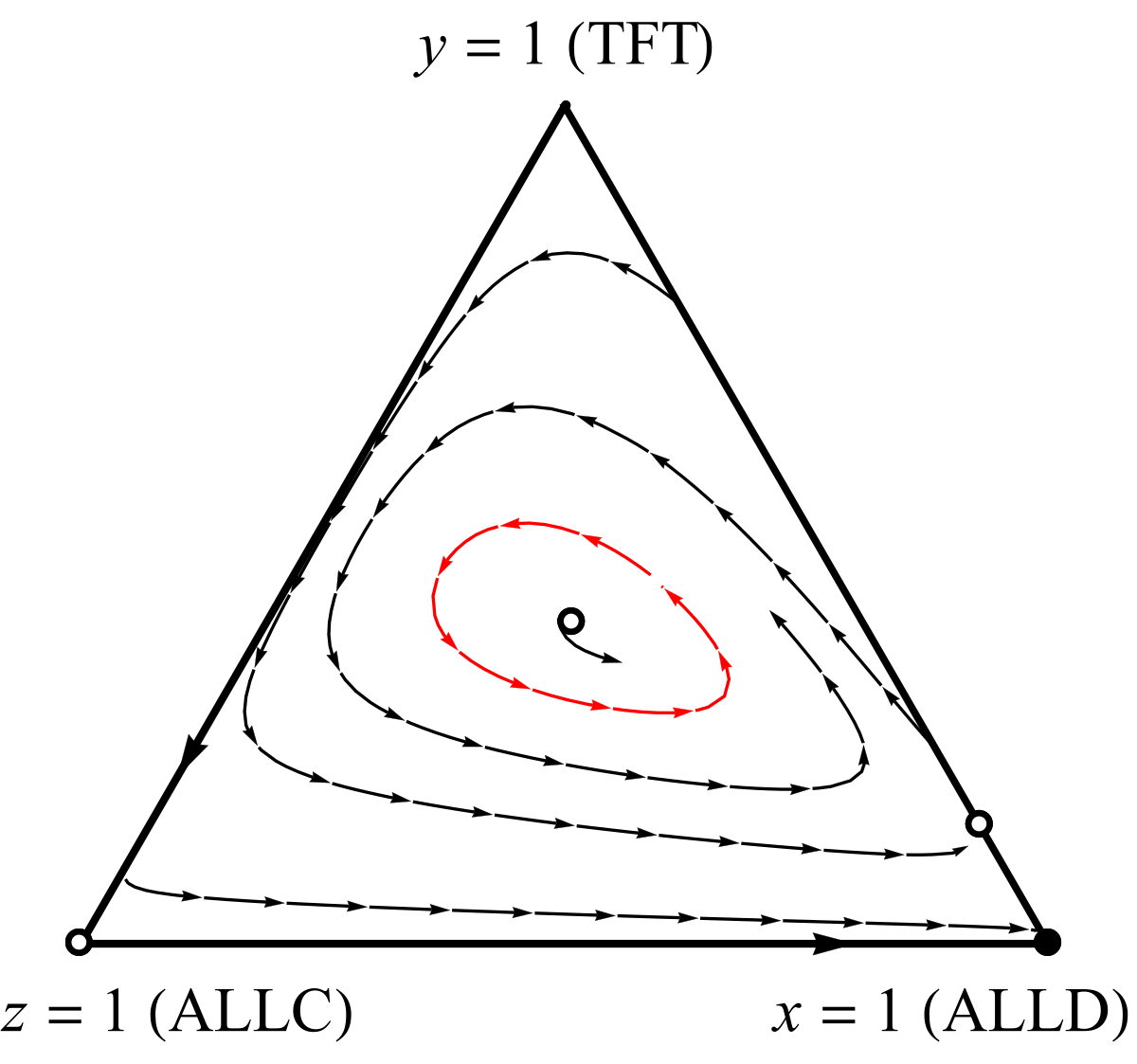}}
%\label{fig:subfigure3}}
\quad
\subfigure[\  Region 4]{%
\includegraphics[scale = 0.13]{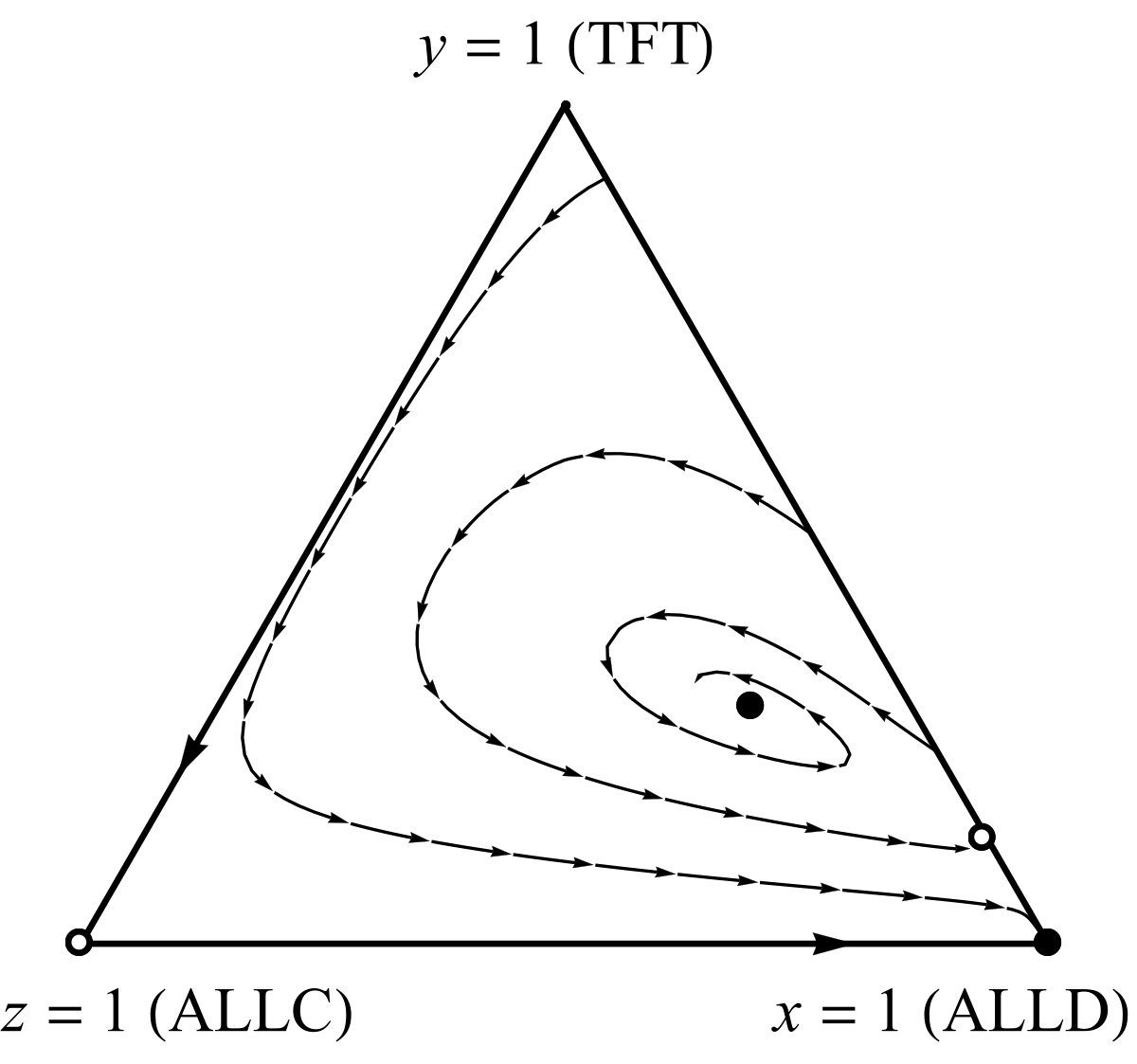}}
%\label{fig:subfigure4}}
%
\caption{Phase portrait of system \eqref{eqn:rmeqn2_TFTALLC} in different regions of $\displaystyle (\mu,c)$ space when TFT $\rightarrow$ ALLC. The stable limit cycle is shown in red.}
\label{fig:TFTALLCphaseportraits}
\end{figure}

In region~3, however, a new attractor -- a stable limit cycle, shown in red in Fig. 4(c) -- coexists with ALLD. Where did this limit cycle come from? It emerged from a homoclinic orbit. When the parameters lie on the homoclinic bifurcation curve between region~2 and region~3, a homoclinic orbit starts and ends at the newly created saddle (the one close to the side between TFT and ALLD). This homoclinic orbit becomes a stable limit cycle when the parameters lie in region~3. 

Finally, as we move from region~3 toward region~4, the limit cycle shrinks and ultimately becomes  a point (a stable spiral) at the supercritical Hopf curve (the red curve in Fig.~\ref{fig:bifnTFTALLC}). When we move into the interior of region~4 the system becomes bistable, with the stable spiral sharing the state space with ALLD. 

In biological terms, the population displays substantial levels of cooperation when it is on either the stable limit cycle of region 3 or the stable spiral of region 4.

%---------  ALLD -> ALLC ------------------------------------------------------------------------------

\subsection{Example 2: ALLD $\rightarrow$ ALLC}
Now we consider an alternative scenario of unidirectional mutation. Suppose that after replication occurs, each player with strategy ALLD mutates into ALLC with probability $\mu$. Notice that this example shares one potentially important feature with Example 1 -- in both cases, the target mutant being created is ALLC. On the other hand, this example differs from Example 1 in that the source of the mutant is now ALLD, not TFT. Which matters more: the commonality of the target or the non-commonality of its source? 

It turns out that the commonality of the target is more important. To see this, observe that in the presence of ALLD $\rightarrow$ ALLC mutations, the replicator-mutator system becomes
\begin{eqnarray}
\dot{x} &=& x \left(f_x - \phi\right) - \mu x f_x, \nonumber \\
\dot{y} &=& y  \left(f_y - \phi\right)    \label{eqn:rmeqn1_ALLDALLC}
\end{eqnarray}
and $z = 1 - x - y$, where $x, y,$ and $z$ again denote the frequencies of ALLD, TFT, and ALLC, respectively.
After insertion of  \eqref{eqn:fitness_with_numbers_and_cost} and  \eqref{eqn:avg_fitness_with_numbers_and_cost} into \eqref{eqn:rmeqn1_ALLDALLC}, the replicator-mutator system yields
\begin{eqnarray}
\dot{x} &=&  x [ \left(c-4\right)  y-5  \mu +4  \mu \left(x+y \right) +x \left(x+3 y-3 \right)+2], \nonumber \\
\dot{y} &=& y  [ c  \left(y-1\right)+x\left(x+3 y-1 \right)].  \label{eqn:rmeqn2_ALLDALLC}
\end{eqnarray}

Figure~\ref{fig:bifnALLDALLC} plots the stability diagram in $(\mu, c)$ space. Note how much it resembles Fig.~\ref{fig:bifnTFTALLC}. In both cases, the parameter space divides into four regions bounded by curves of supercritical Hopf (red), saddle-node (blue) and homoclinic (green) bifurcations, all of which emerge from a Bogdanov-Takens point. The main difference is that the red curve of Hopf bifurcations goes through the origin in Fig.~\ref{fig:bifnTFTALLC} and not in Fig.~\ref{fig:bifnALLDALLC}.

\begin{figure}[h]
\centering
\includegraphics[scale = 0.27]{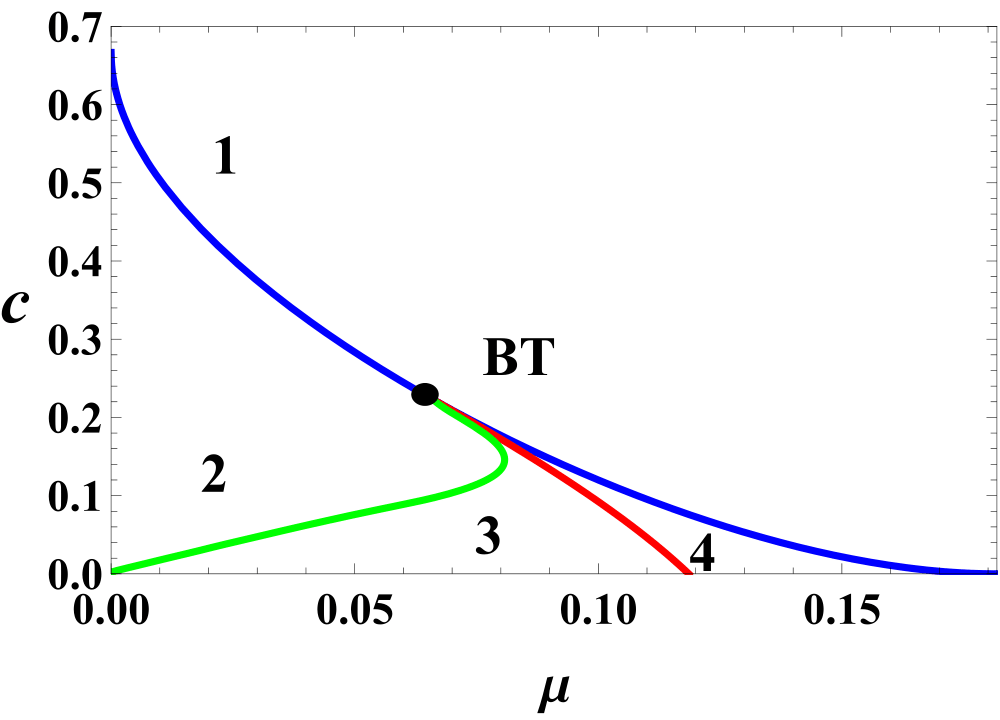}
\caption{Stability diagram in $(\mu, c)$ space when ALLD $\rightarrow$ ALLC. Saddle-node curve, blue; supercritical Hopf curve, red; homoclinic curve, green; BT, Bogdanov-Takens point. Formulas for the Hopf and saddle-node bifurcation curves are given in  Appendix B. The homoclinic curve was computed using MATCONT.} 
\label{fig:bifnALLDALLC}
\end{figure}

Figure~\ref{fig:ALLDALLCphaseportraits} shows the phase portraits corresponding to the four regions. As in Example~1, defection  dominates the long-term dynamics in regions 1 and 2, where a stable fixed point close to ALLD attracts almost all solutions. (Note that pure ALLD is no longer a fixed point, because of the assumed mutations from ALLD  to ALLC. That is why the globally attracting fixed point lies between ALLD and ALLC. We will refer to this fixed point as ``almost ALLD.'') In region~3, a stable limit cycle (shown in red) coexists with almost ALLD. In region~4, the limit cycle no longer exists; it has contracted to a stable spiral fixed point via the supercritical Hopf bifurcation between regions~3 and 4.

\begin{figure}[h!]
\centering
\subfigure[\  Region 1]{%
\includegraphics[scale = 0.17]{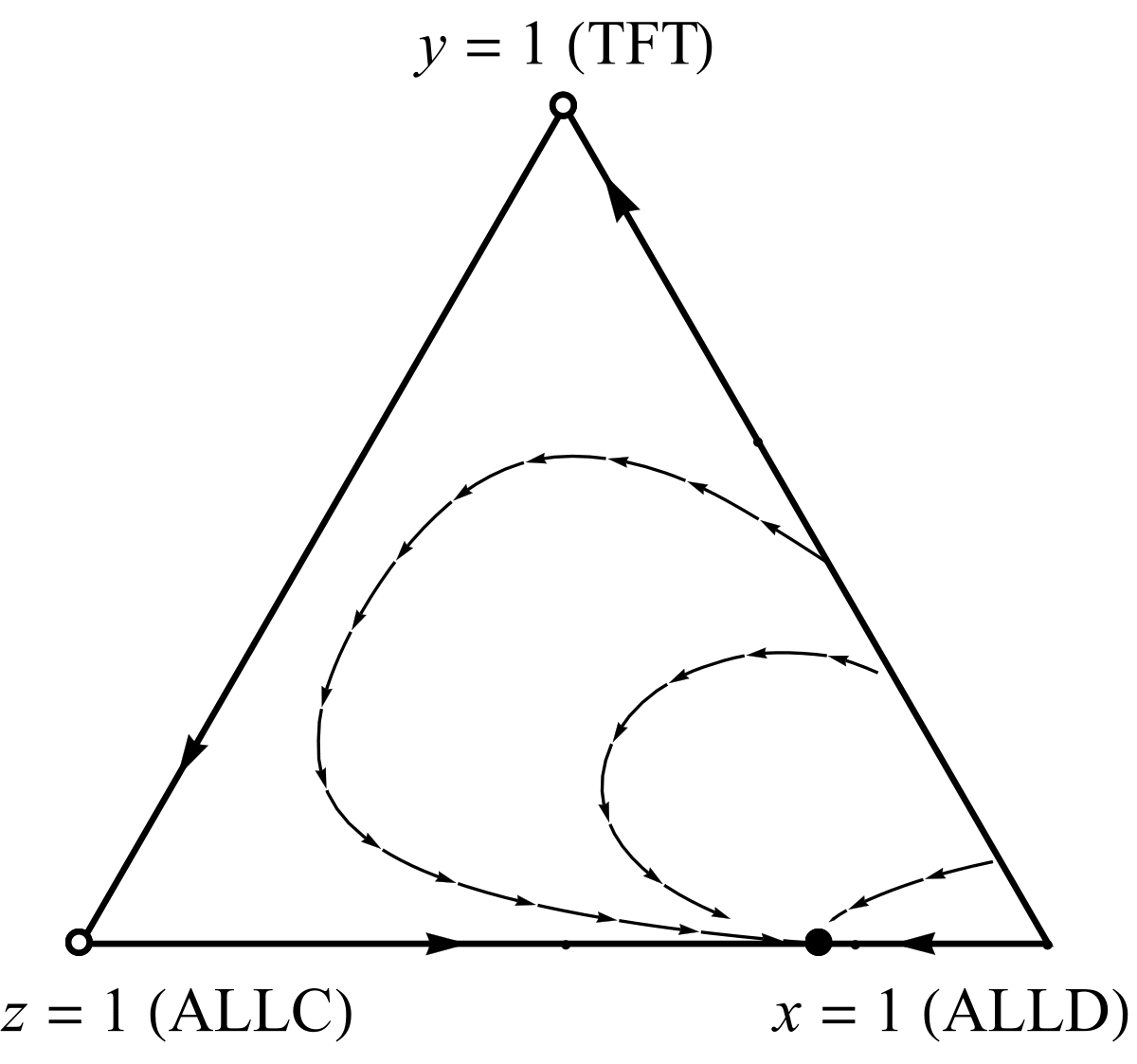}}
%\label{fig:subfigure1}}
\quad
\subfigure[\  Region 2]{%
\includegraphics[scale = 0.17]{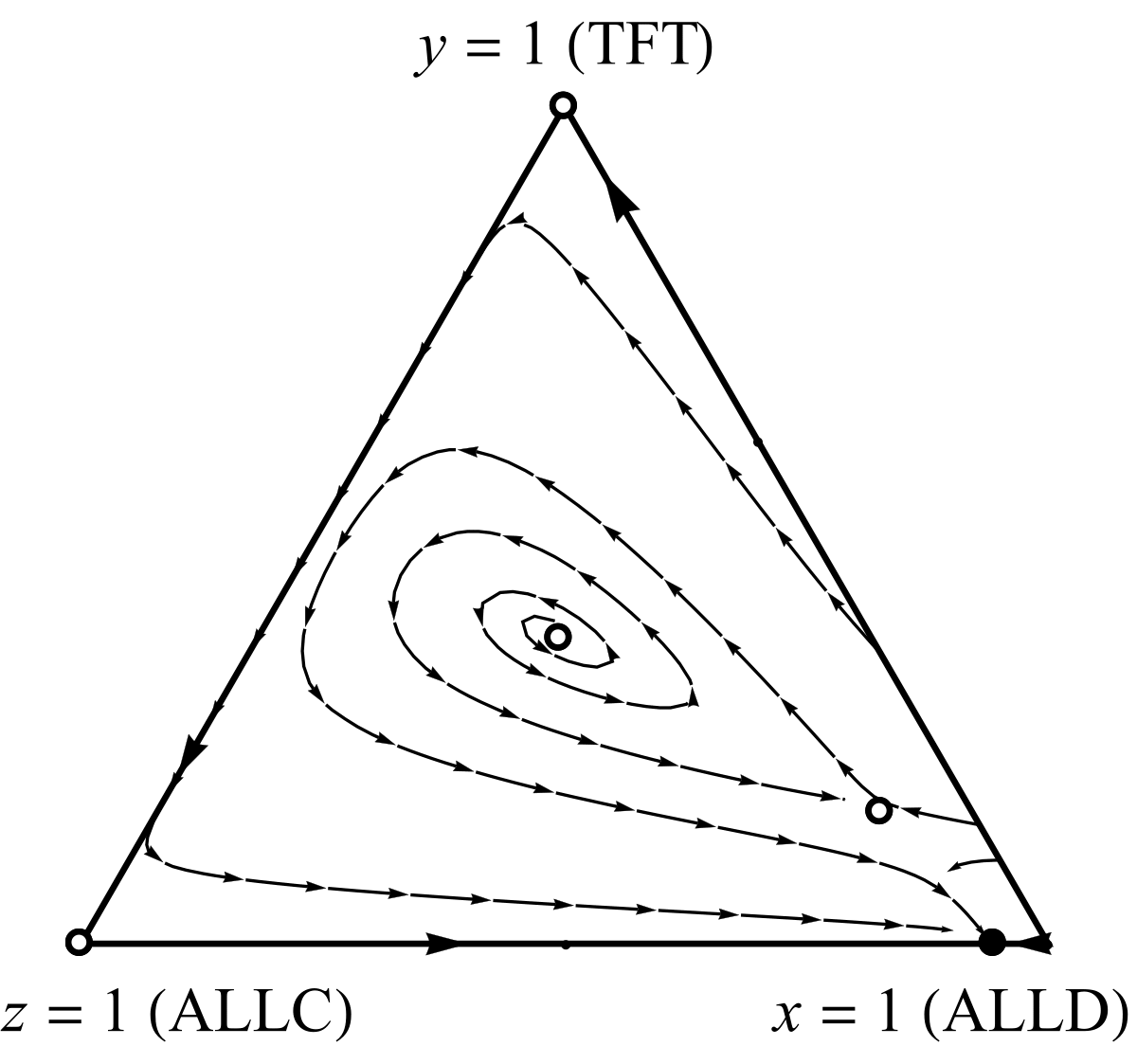}}
%\label{fig:subfigure2}}

\subfigure[\  Region 3]{%
\includegraphics[scale = 0.2]{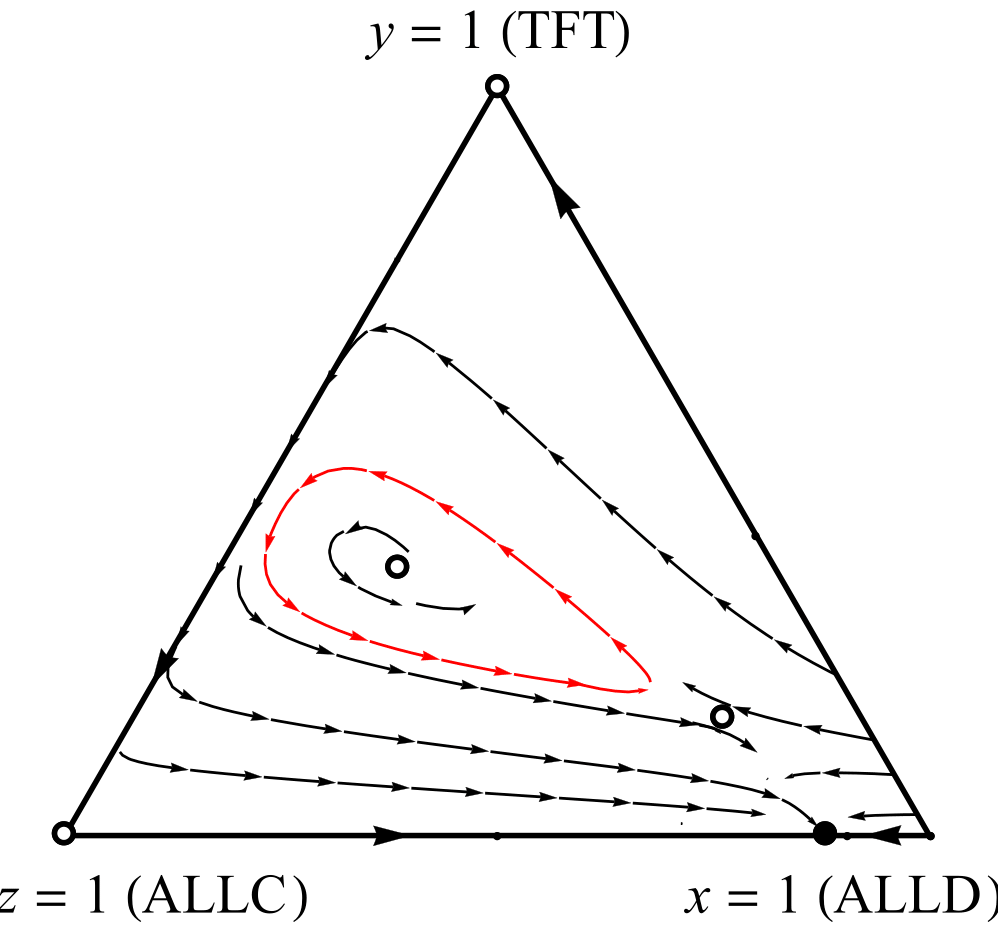}}
%\label{fig:subfigure3}}
\quad
\subfigure[\  Region 4]{%
\includegraphics[scale = 0.17]{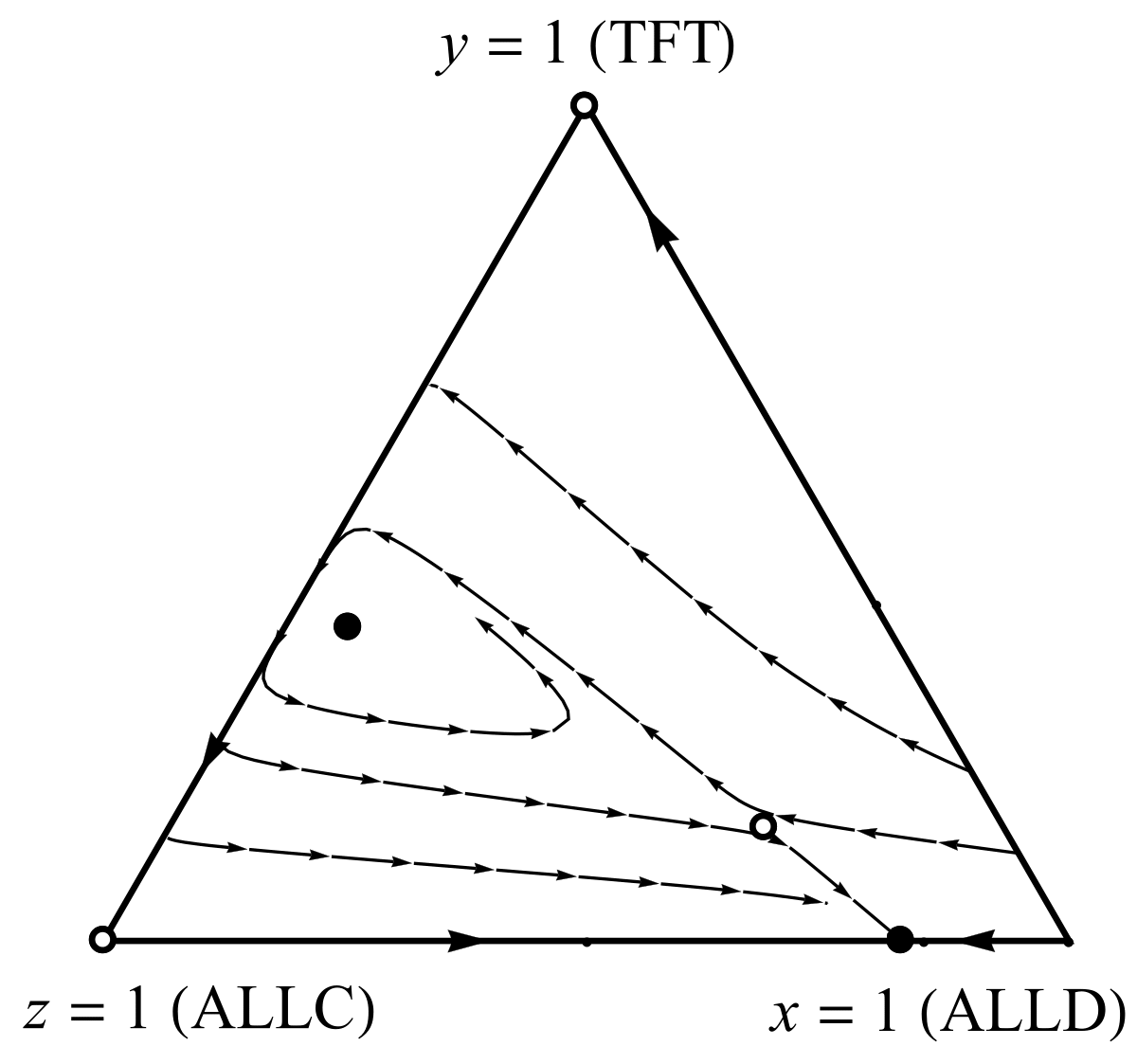}}
%\label{fig:subfigure4}}
%
\caption{Phase portrait of system \eqref{eqn:rmeqn2_ALLDALLC} in different regions of $(\mu, c)$ space when ALLD $\rightarrow$ ALLC. The stable limit cycle is shown in red.}
\label{fig:ALLDALLCphaseportraits}
\end{figure}

It is important to realize that the population exhibits large amounts of cooperation when it is in the stable spiral state or cycling periodically. This becomes clear when one looks at time series instead of phase portraits. Figure~\ref{fig:ALLDALLCtimeseries} shows the approach to a stable limit cycle in which  ALLD is nearly absent except during brief spikes, whereas ALLC and TFT predominate for most of the cycle. Thus, the stable limit cycle signifies more than just an avoidance of all-out defection; it represents a state of significant cooperation.  

\begin{figure}[h!]
\centering
\subfigure[\  ALLD]{%
\includegraphics[scale = 0.17]{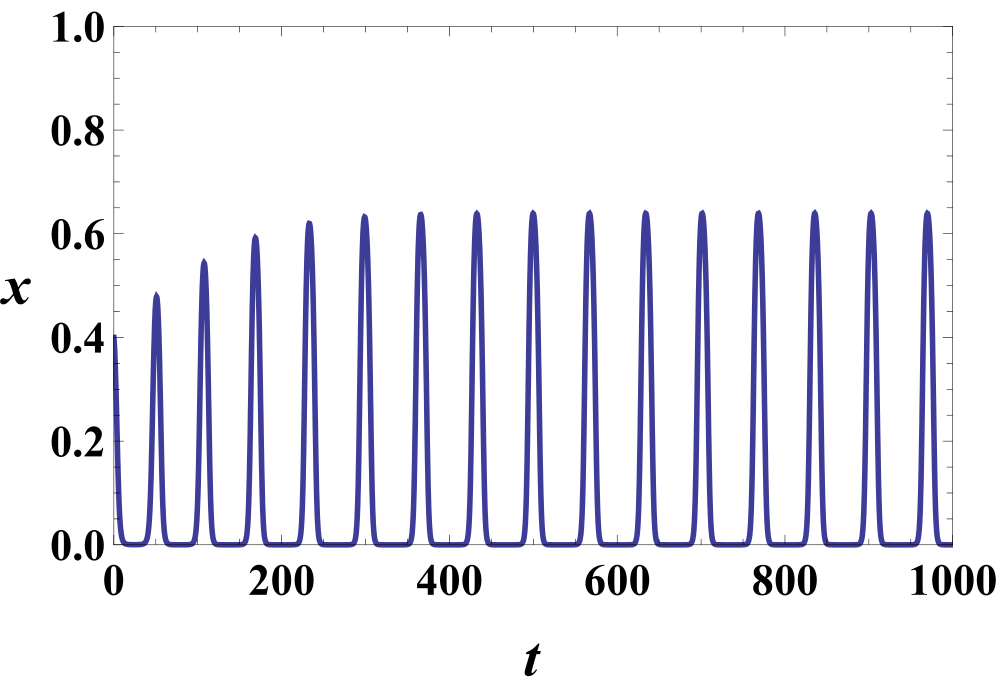}}
%\label{fig:subfigure1}}
\subfigure[\  TFT]{%
\includegraphics[scale = 0.17]{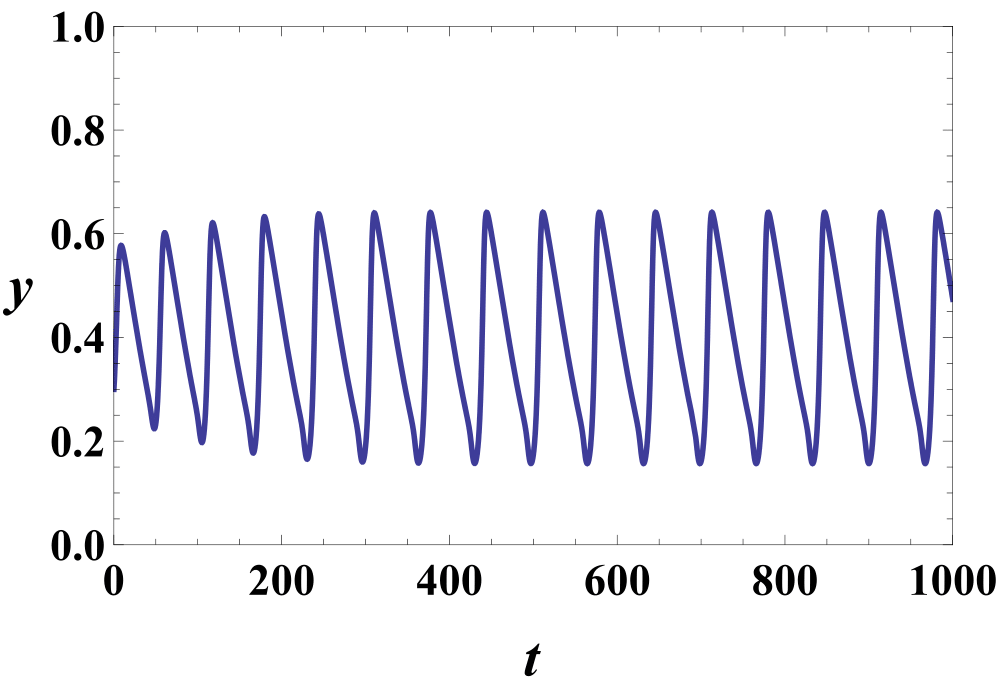}}
%\label{fig:subfigure2}}
\subfigure[\  ALLC]{%
\includegraphics[scale = 0.17]{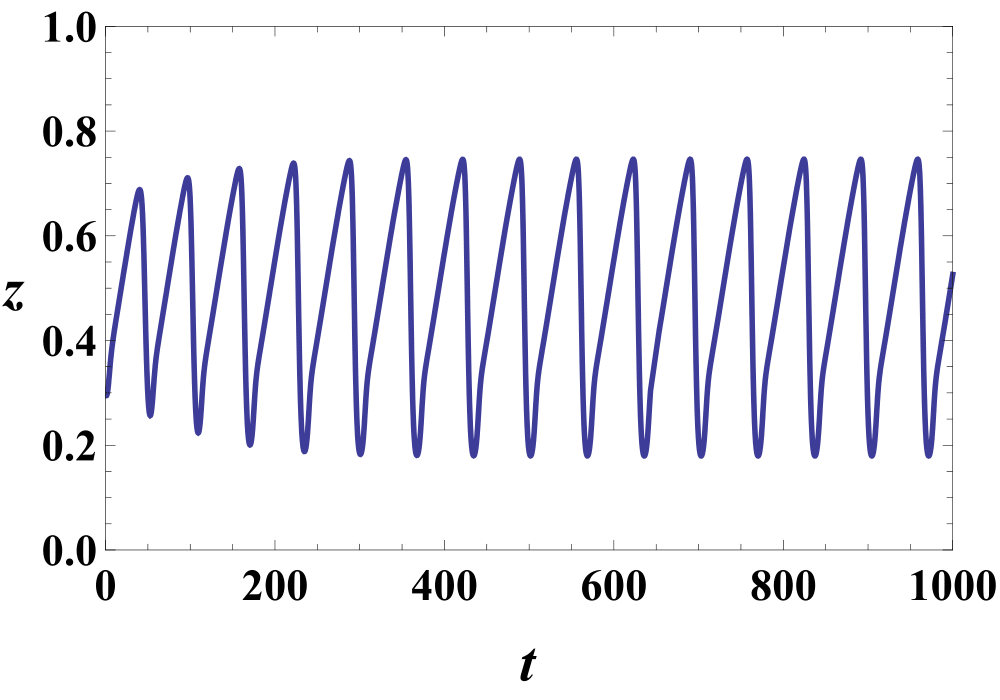}}
%\label{fig:subfigure3}}
\caption{Time series for a solution of system \eqref{eqn:rmeqn2_ALLDALLC} as it approaches the stable limit cycle. Notice that the level of cooperation is high during much of the cycle, as shown by the high combined levels of ALLC and TFT, while ALLD remains low except for brief spikes. Parameter values: $(\mu, c) = (0.08, 0.04)$.}
\label{fig:ALLDALLCtimeseries}
\end{figure}

%------------------------------  Other single unidirectional mutations------------------
\subsection{Other single unidirectional mutations}
Four other types of single unidirectional mutations are possible: ALLD~$\rightarrow$ TFT, ALLC~$\rightarrow$ TFT, ALLC~$\rightarrow$ ALLD, and TFT~$\rightarrow$ ALLD. Using the techniques above, we have analyzed these remaining cases completely. But rather than wade through the details, it seems clearer and more useful to summarize the main results, which are as follows. 

For all four cases, the phase portraits and the bifurcation curves are more complicated than those presented in Examples 1 and 2. Nevertheless, all of them display stable limit cycles in some region of the parameter space. Specifically, the region in $(\mu, c)$ space where stable limit cycles exist is always bounded on one side by a curve of supercritical Hopf bifurcations, and on the other side by a curve of homoclinic bifurcations. Both  curves emanate from a Bogdanov-Takens point. 

What we find most intriguing is that the homoclinic bifurcation curve -- the counterpart of the green curve in Figs.~\ref{fig:bifnTFTALLC} and \ref{fig:bifnALLDALLC} -- always passes through the origin $(\mu, c) = (0,0)$. Hence stable limit cycles always exist for arbitrarily small perturbations of the original system~\eqref{eqn:replicator} with fitnesses \eqref{eqn:fitness_with_numbers}, no matter how the mutation is implemented. 

This is the sense in which limit cycles are ``sparked by mutation'' in the Prisoner's Dilemma among ALLC, ALLD, and TFT. If just one of the strategies can mutate into just one of the others,  it takes only an infinitesimal cost $c$ and an infinitesimal mutation probability $\mu$ for the system to display self-sustained oscillations, with large amounts of cooperation during part of the cycle. In this way, the slightest bit of mutation allows the population to avoid a collapse into all-out defection.

%------------------------------  Global------------------------------------------------------------------------------

\subsection{Uniform global mutation}
So far we have focused on a very restricted class of mutation pathways: single unidirectional mutations, in which exactly one strategy mutates into exactly one other. But we suspect that the sparking of limit cycles by mutation is more general. 

For example, consider the extreme opposite case of uniform global mutation, where each strategy mutates to the other two with probability $\mu$, and hence stays the same with probability $1-2 \mu$. The replicator-mutator equations for this case are given by
\begin{eqnarray}
\dot{x} &=& x \left[f_x \left(1-2   \mu\right) - \phi \right] +\mu   f_y   y+\mu   f_z   z \nonumber\\
\dot{y} &=& y  \left[f_y  \left(1-2   \mu\right) - \phi \right]  +\mu   f_x x+\mu   f_z  z  \label{eqn:rmeqn1_global}
\end{eqnarray}
where $z=1-x-y$. 
Substitution of \eqref{eqn:fitness_with_numbers_and_cost} and  \eqref{eqn:avg_fitness_with_numbers_and_cost} into \eqref{eqn:rmeqn1_global} yields
\begin{eqnarray}
\dot{x} &=& \mu  \left[x \left(11  x+9  y-16\right)-c  y \right] +x \left[\left(c-4\right)  y+x\left(x+3  y-3\right)+2 \right]+3 \mu, \nonumber\\
\dot{y} &=& c  y \left(2  \mu +y-1\right)+3  \mu +x^2  \left(y-\mu \right) +x \left(3  y-1\right)  \left(\mu +y\right)-9   \mu   y.
 \label{eqn:rmeqn2_global} 
\end{eqnarray} 

Figure~\ref{fig:globalbifurcation} plots the stability diagram in $(\mu, c)$ space, showing the curves where supercritical Hopf (red), saddle-node (blue) and homoclinic (green) bifurcations take place. The diagram splits into five regions. All the bifurcation curves of this system were computed using MATCONT [Govaerts \&  Kuznetsov, 2008].

\begin{figure}[h]
\centering
\includegraphics[scale = 0.27]{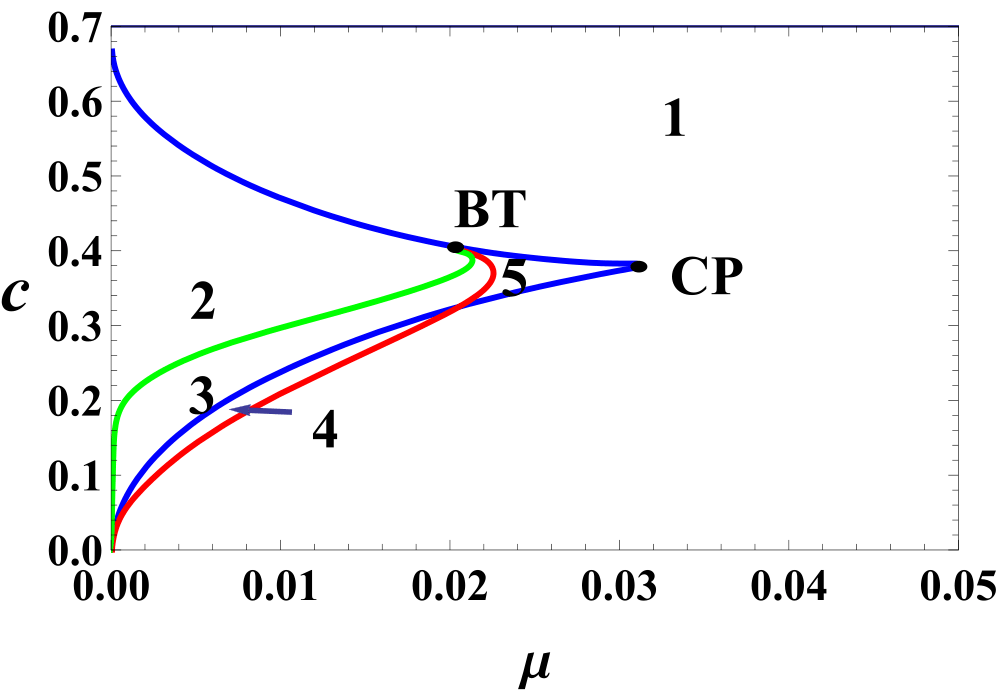}
\caption{Stability diagram in ($\mu, c)$ space for the case of uniform global mutation. Saddle-node curve, blue; supercritical Hopf curve, red; homoclinic curve, green; BT, Bogdanov-Takens point; CP, cusp point. The bifurcation curves were computed using MATCONT.}  
\label{fig:globalbifurcation}
\end{figure}

Figure~\ref{fig:globalphaseportraits} shows the phase portraits corresponding to the five regions. As in the earlier examples~1 and 2, defection prevails in regions 1 and 2, where a stable fixed point close to ALLD attracts almost all solutions.  In region~3, the possibility of cooperation reemerges: a stable limit cycle (shown in red) coexists with a stable node near ALLD. In region~4, the limit cycle becomes a global attractor; the stable node of region~3 is lost in a saddle-node bifurcation when the parameters move from region~3 to 4.  The limit cycle no longer exists in region~5. 

\begin{figure}[!h]
\centering
\subfigure[\  Region 1]{%
\includegraphics[scale = 0.17]{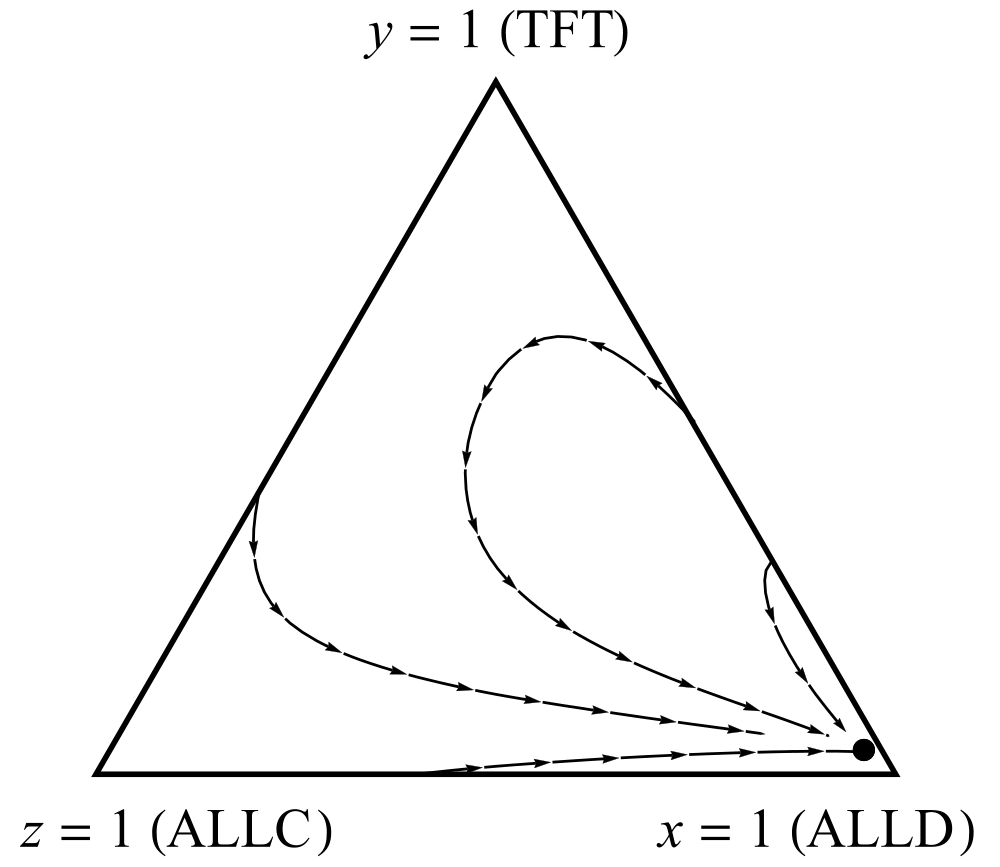}}
%\label{fig:subfigure1}}
\quad
\subfigure[\  Region 2]{%
\includegraphics[scale = 0.17]{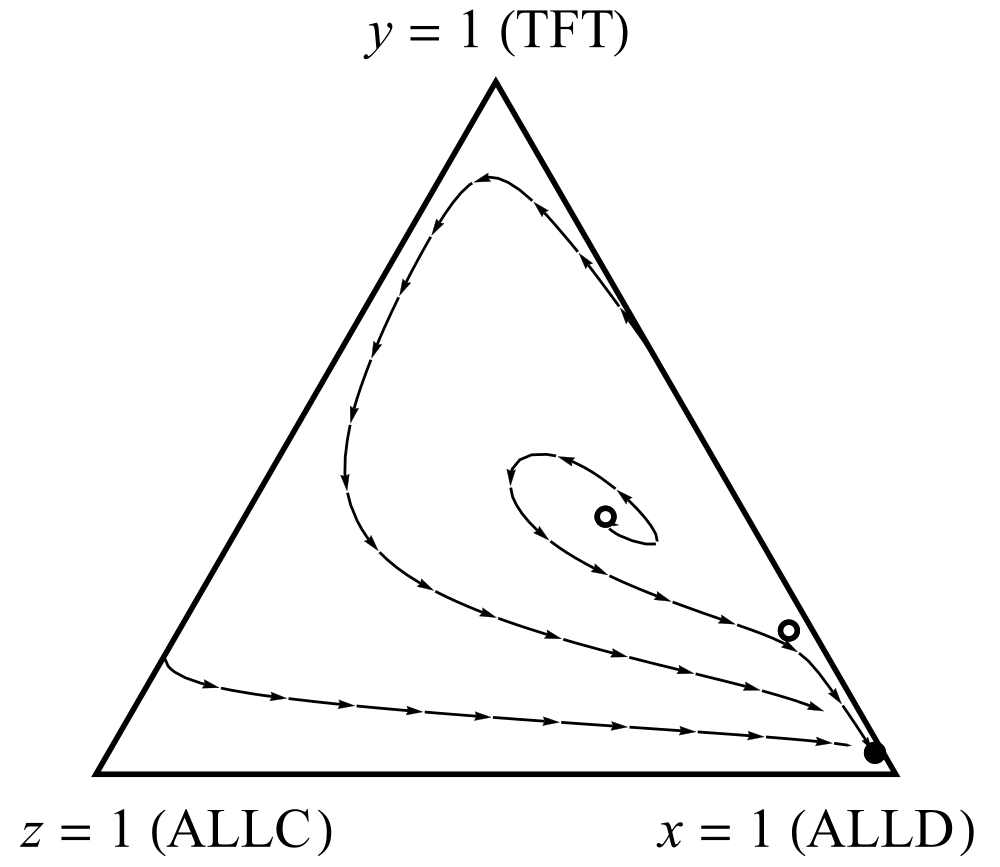}}
%\label{fig:subfigure2}}

\subfigure[\  Region 3]{%
\includegraphics[scale = 0.17]{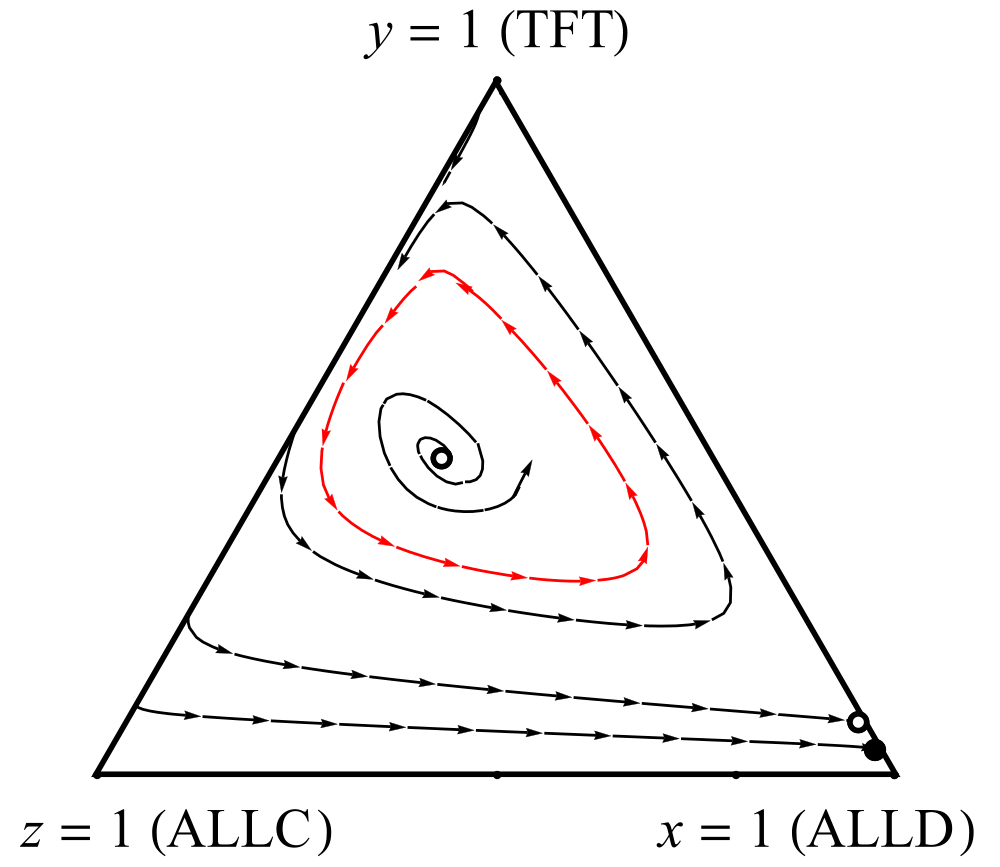}}
%\label{fig:subfigure3}}
\quad
\subfigure[\  Region 4]{%
\includegraphics[scale = 0.17]{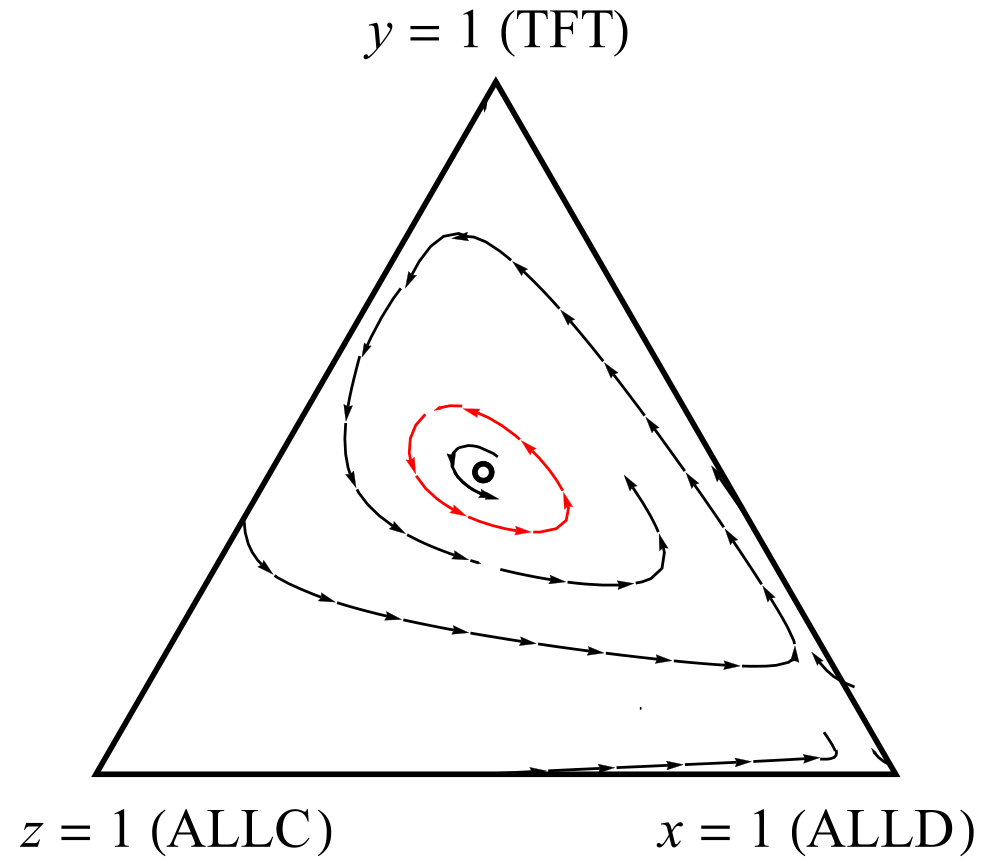}}
%\label{fig:subfigure4}}
\quad
\subfigure[\  Region 5]{%
\includegraphics[scale = 0.17]{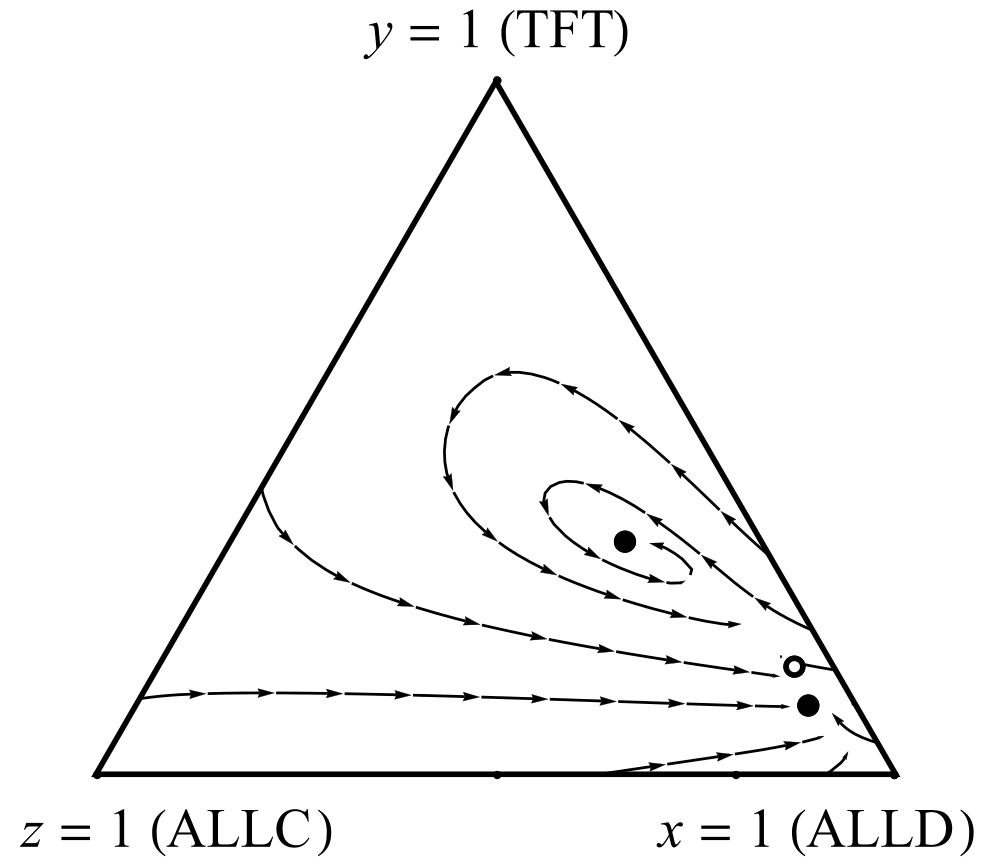}}
%\label{fig:subfigure4}}
%
\caption{Phase portrait of system \eqref{eqn:rmeqn2_global} in different regions of $(\mu,c)$ parameter space. The stable limit cycle is shown in red.}
\label{fig:globalphaseportraits}
\end{figure}

\subsection{Conjecture}
In every one of the examples considered so far, the region where stable limit cycles exist extends all the way down to the origin in parameter space. This means that limit cycles can be sparked by an arbitrarily small mutation probability $\mu$, if the complexity cost $c$ of TFT is also suitably small. These results have been obtained for topologies at opposite ends of the graph-theoretic spectrum: single unidirectional mutation (in which mutation occurs along one directed edge) and uniform global mutation (in which mutation occurs in both directions along the edges of the complete graph). 

We conjecture, therefore, that a similar sparking of stable limit cycles occurs for \emph{any} pattern of mutation. To make this statement precise, we write down the general form of the replicator-mutator equation. Let the frequencies of the three strategies ALLD, TFT and ALLC be denoted as $x_1, x_2, x_3$ rather than $x, y,z$. Then the replicator-mutator dynamics are given by

\begin{eqnarray}
\dot{x_i} &=& \left( \sum_{j=1}^3  x_j f_j Q_{ji} \right) - x_i \phi \label{eqn:rm_general}
\end{eqnarray}

\noindent for $i = 1, 2, 3$. Here $f_i$ is the fitness of strategy $i$ (obtained from Eq.~\eqref{eqn:fitness_with_numbers}, as before), $\phi = \sum_{j=1}^3  x_j f_j$ is the average fitness in the population, and $Q_{ij}$ is the probability that players with strategy $i$ mutate to playing strategy $j$. Since the $Q_{ij}$ represent probabilities, they are non-negative and satisfy $\sum_j Q_{ij} = 1$; hence the matrix $Q$ is row stochastic. In the limiting case where mutation does not occur, $Q$ is the identity matrix and $Q_{ij} = \delta_{ij}$ where $\delta_{ij}$ is the Kronecker delta. 

Now suppose that mutation does occur, and that it can be characterized by a mutation matrix $M$ and a single parameter $0 < \mu \leq 1$ such that 
\begin{eqnarray}
Q_{ij} = \delta_{ij} - \mu M_{ij}. \label{eqn:M}
\end{eqnarray}
Any matrix $M$ that maintains the row-stochasticity and non-negativity of $Q$ is admissible. This diversity of $M$ is what we mean by any pattern of mutation. 

Phrased in these terms, our conjecture is that for any fixed, admissible $M$, one can find an open set of parameters $(\mu,c)$ arbitrarily close to $(0,0)$ such that stable limit cycles exist for the replicator-mutator system \eqref{eqn:rm_general}.    

We have verified this conjecture informally for a handful of admissible mutation matrices $M$, but we have not done a comprehensive numerical study nor do we have analytical evidence to support the conjecture. So it could well be false. But we suspect it is true. A key step toward proving it would be to demonstrate that a curve of homoclinic bifurcations (the counterpart of the green curves in Figs.~\ref{fig:bifnTFTALLC}, \ref{fig:bifnALLDALLC}, and \ref{fig:globalbifurcation}) always emanates from the point $(\mu,c) = (0,0)$, for any admissible $M$.

%-----------------------------------Discussion----------------------------------------------------------------------------

\section{Discussion}
Our results show that a variety of different mutation structures give rise to evolutionary cycles of cooperation and defection in the repeated Prisoner's Dilemma. These cycles appear to be a robust feature of interactions in well-mixed populations of ALLC, ALLD and TFT, rather than being specific to particular assumptions regarding the mutation matrix. These cycles, as well as the stable spirals created by mutation in other parts of the $(\mu, c)$ parameter space, lead to substantial levels of cooperation. Thus the grim picture painted by the standard replicator equation without mutation, in which ALLD dominates even in repeated games, may be too pessimistic. 

The replicator equation assumes that the population is well mixed. In some settings it is more realistic to regard the population as spatially structured, with players interacting only with a fixed subset of others on a regular lattice. Evolutionary cycles of cooperation and defection among ALLC, ALLD, and TFT can occur in this case too, even in the absence of mutation \cite{Szolnoki09}.  Similar cycles occur for populations playing a spatial version of the rock-paper-scissors game  \cite{Szolnoki04}.

Although the replicator equation is typically thought of as describing genetic evolution, it can just as well describe a social learning dynamic in which people imitate the strategies of more successful others. In this context, mutation corresponds to experimentation or innovation, i.e., trying out a strategy other than the one which is performing well at the moment. Experimentation is a key element of human behavior \cite{Rand132, Traulsen10}, and our results suggest that it may also help to explain the success of cooperation in human societies.

%-----------------------------------Acknowledgments----------------------------------------------------------------------------

\section{Acknowledgments} Research supported in part by a Sloan-Colman fellowship to D.F.P.T. and a grant from the John Templeton Foundation through a subaward from the New Paths to Purpose Project at Chicago Booth to D.G.R. We thank Richard Rand and Elizabeth Wesson for helpful discussions.

\clearpage
%------------------------References--------------------------------------------
%\nonumsection{References}

%
%\end{multicols}

%--------------------------APPENDICES-------------------------------------------------------------------------------------------
\clearpage
\appendix{TFT $\rightarrow$ ALLC}

For Example 1 of the main text, in which TFT mutates into ALLC, the system \eqref{eqn:rmeqn1_TFTALLC} has at most 4 fixed points $(x^*,y^*)$ in the simplex $x+y+z = 1$, with $0 \leq x,y,z \leq 1$:\\
%\begin{equation*}
%\left(
%\begin{array}{cc}
% \displaystyle x = 0 &  \displaystyle y= \left(\frac{3}{c}-1\right) \mu +1 \\
% \displaystyle x = 1 &  \displaystyle y = 0 \\
% \displaystyle x = \frac{-A_1-5 c \mu +c+17 \mu +2}{12 \mu +4} &  \displaystyle y = \frac{\mu  \left(A_1+\mu
%   -1\right)+A_1+c (1-(\mu -8) \mu )+2}{24 \mu +8} \\
% \displaystyle x = \frac{A_1-5 c \mu +c+17 \mu +2}{12 \mu +4} &  \displaystyle y = \frac{\mu  \left(-A_1+\mu
%   -1\right)-A_1+c (1-(\mu -8) \mu )+2}{24 \mu +8} \\
% \displaystyle x = 0 &  \displaystyle y = 0 
%\end{array}
%\right)
%\end{equation*} \\
%\begin{equation*}
%\begin{align}

$ \displaystyle (x^*_1,y^*_1) = \displaystyle \left(0,0\right)$  is a saddle, \\

$ \displaystyle (x^*_2,y^*_2) =  \displaystyle (1,0)$  is a stable node,\\

$ \displaystyle (x^*_3,y^*_3) =  \displaystyle \left(\frac{A_1-5 c \mu +c+17 \mu +2}{12 \mu +4},\frac{-(\mu +1) A_1-c \mu ^2+8 c \mu +c+\mu ^2-\mu +2}{24 \mu +8}\right)$  is a saddle, and\\

$ \displaystyle (x^*_4,y^*_4) = \displaystyle \left(\frac{-A_1-5 c \mu +c+17 \mu +2}{12 \mu +4},\frac{(\mu +1) A_1-c \mu ^2+8 c \mu +c+\mu ^2-\mu +2}{24 \mu +8}\right)$,  could be a stable or unstable spiral (or node), depending on  $\mu$ \text{and} $c$ \cite{Strogatz94}, where in the formulas above,\\

$\displaystyle A_1 = \sqrt{c^2 (\mu +3)^2-2 c ((\mu -11) \mu +6)+(\mu -28) \mu +4}.$\\

The fixed points  $(x^*_1,y^*_1)$, $(x^*_2,y^*_2) $, and $(x^*_3,y^*_3)$  are on the boundary of the simplex and $(x^*_4,y^*_4)$ is inside the simplex. A saddle-node bifurcation occurs when $(x^*_3,y^*_3)$ and $(x^*_4,y^*_4)$ coalesce and disappear. The equation of the saddle-node bifurcation curve is 
\begin{equation}
c = \displaystyle \frac{(\mu -11) \mu -4 \sqrt{6} \sqrt{\mu  (3 \mu +1)}+6}{(\mu +3)^2}. \label{eqn:TFTALLC_SN}
\end{equation}

The fixed point $(x^*_4,y^*_4)$ undergoes a supercritical Hopf bifurcation and switches from a stable spiral to an unstable spiral at certain parameter values. The equation of the Hopf bifurcation curve, which was computed analytically, is \\
%\begin{equation}
\begin{align}
c = \displaystyle -\frac{\displaystyle 3-13 \mu }{\displaystyle 4 \left( \mu -1 \right)} &-\frac{1}{2}\ \sqrt{\displaystyle A_2+\frac{A_4}{\displaystyle 3 \sqrt[3]{2} (\mu
   -1)^2}+\frac{A_5}{\displaystyle 3 (\mu -1)^2 A_4}} \\ &+  \frac{1}{2} \sqrt{\displaystyle A_2-\frac{A_4}{\displaystyle 3
   \sqrt[3]{2} (\mu -1)^2}-\frac{A_5}{\displaystyle 3 (\mu -1)^2 A_4}+ \frac{A_6}{\displaystyle 4 \sqrt{A_7+\frac{A_4}{\displaystyle 3 \sqrt[3]{2} (\mu
   -1)^2}+\frac{A_5}{\displaystyle 3 (\mu -1)^2 A_4}}}} \label{eqn:TFTALLC_Hopf}
\end{align}
%\end{equation}

\noindent where\\

\smallskip

$\displaystyle A_2 = \frac{(3-13 \mu )^2}{4 (\mu -1)^2}+ \frac{35 \mu ^2-34 \mu -8}{3 \left(\mu ^2-2 \mu +1\right)}-\frac{35 \mu ^2-34 \mu -8}{(\mu
   -1)^2},$\\

\smallskip

$\displaystyle A_3 = (-5240604096 \mu ^{12}-40578465024 \mu ^{11}+200756188800 \mu ^{10}-265354820640 \mu
   ^9+60401533248 \mu ^8+110419920576 \mu ^7-55059298560 \mu ^6-12327872544 \mu ^5+4987630080 \mu ^4+1779338880 \mu ^3+218439936 \mu ^2-110592 \mu -1880064)^{\frac{1}{2}},$\\

\smallskip

$\displaystyle A_4 = (204136 \mu ^6-578832 \mu ^5+416952 \mu ^4+51238 \mu ^3-79740 \mu ^2-12984 \mu +A_3-1456)^{\frac{1}{3}},$\\

\smallskip

$\displaystyle A_5 =  \sqrt[3]{2} \left(2272 \mu ^4-3160 \mu ^3+521 \mu ^2+316 \mu +100\right),$\\

\smallskip

$\displaystyle A_6 =  -\frac{(3-13 \mu )^3}{(\mu -1)^3}+\frac{4 \left(35 \mu ^2-34 \mu -8\right) (3-13 \mu )}{(\mu -1)^3}-\frac{8 \left(49 \mu ^2-4 \mu +4\right)}{(\mu -1)^2},$ and \\

\smallskip

$\displaystyle A_7 = \frac{35 \mu ^2-34 \mu -8}{3 \left(\mu ^2-2 \mu +1\right)}-\frac{35 \mu ^2-34 \mu -8}{(\mu
   -1)^2}.$

%%%%%% APPENDIX  B %%%%%%%%

\clearpage
\appendix{ALLD $\rightarrow$ ALLC}
For Example 2 of the main text, in which ALLD mutates into ALLC, system \eqref{eqn:rmeqn1_ALLDALLC} has at most 5 fixed points $(x^*,y^*)$ in the simplex $x+y+z = 1$, with $0 \leq x,y,z \leq 1$:\\
%\begin{equation*}
%\left(
%\begin{array}{cc}
% x = \displaystyle  0 & y =\displaystyle 1 \\
% x =\displaystyle \frac{1}{2} \left(-4 \mu -\sqrt{4 \mu  (4 \mu -1)+1}+3\right) & y =\displaystyle 0 \\
% x = \displaystyle -\frac{4 c \mu +A_8+c-11 \mu +2}{16 \mu -4} & y =\displaystyle \frac{8 c \mu ^2-10 c \mu +2 \mu A_8-A_8+c+18 \mu ^2-11 \mu +2}{8 (\mu -1) (4 \mu
%   -1)} \\
% x =\displaystyle \frac{-4 c \mu +A_8-c+11 \mu -2}{16 \mu -4} & y =\displaystyle \frac{8 c \mu ^2-10 c \mu -2 \mu  A_8+A_8+c+18 \mu ^2-11 \mu +2}{8 (\mu -1) (4 \mu
%   -1)} \\
% x =\displaystyle 0 & y =\displaystyle 0 \\
%\end{array}
%\right)
%\end{equation*}
%\begin{equation*}
%\begin{align}

$ \displaystyle(x^*_1,y^*_1) = \displaystyle \left(0,0\right)$   is a saddle, \\

 $\displaystyle(x^*_2,y^*_2) =  \displaystyle (0,1)$  is a saddle,\\

 $\displaystyle(x^*_3,y^*_3)  = \displaystyle \left(\frac{1}{2} \left(-4 \mu -\sqrt{4 \mu  (4 \mu -1)+1}+3\right),0 \right)$  is a stable node,

 $\displaystyle(x^*_4,y^*_4) = \displaystyle \left(  -\frac{4 c \mu +A_8+c-11 \mu +2}{16 \mu -4}  ,\frac{8 c \mu ^2-10 c \mu +(2 \mu -1) A_8+c+18 \mu ^2-11 \mu +2}{8 (\mu -1) (4 \mu -1)}\right)$  is a saddle, and\\

$ \displaystyle (x^*_5,y^*_5) = \displaystyle \left( \frac{-4 c \mu +A_8-c+11 \mu -2}{16 \mu - 4} ,  \frac{8 c \mu ^2-10 c \mu -(2 \mu -1) A_8+c+18 \mu ^2-11 \mu +2}{8 (\mu -1) (4 \mu -1)}\right)$ could be a stable or unstable spiral (or node), depending on   $\mu$ \text{and} $c$. In the expressions above,\\

$\displaystyle A_8 = \sqrt{(4 c \mu +c-11 \mu +2)^2+8 c (4 \mu -1) (-c+\mu +2)}.$\\

 The fixed points  $(x^*_1,y^*_1)$, $(x^*_2,y^*_2) $, and $(x^*_3,y^*_3)$  are on the boundary of the simplex, and $(x^*_4,y^*_4)$ and $(x^*_5,y^*_5)$  are in the interior. A saddle-node bifurcation occurs when $(x^*_4,y^*_4)$ and $(x^*_5,y^*_5)$ coalesce and disappear. The equation of the saddle-node bifurcation curve is

\begin{equation}
c =  \frac{\mu\  (28 \mu -25)+6-4 \sqrt{6}\sqrt{(1-\mu) \mu  (3 \mu -2) (4 \mu -1)}}{(3-4 \mu )^2}. \label{eqn:ALLDALLC_SN}
\end{equation}

The fixed point $(x^*_5,y^*_5)$ undergoes a supercritical Hopf bifurcation and switches from a stable spiral to an unstable spiral at the curve given by \\
\begin{equation}
c = \displaystyle-\frac{\left(1-i \sqrt{3}\right) A_{10}}{6 \sqrt[3]{2}}+\frac{\left(1+i \sqrt{3}\right) \left(48 \mu ^3-201 \mu ^2+114 \mu-33\right)}{3\ 2^{2/3} A_{10}}+2 \mu +1  \label{eqn:ALLDALLC_Hopf}
\end{equation}

\noindent where\\

$\displaystyle A_9 = \sqrt{4 \left(48 \mu ^3-201 \mu ^2+114 \mu -33\right)^3+\left(-2592 \mu ^4+5670
   \mu ^3-5022 \mu ^2+810 \mu +162\right)^2}$\\

\smallskip 

\noindent and
\smallskip 

$\displaystyle A_{10}= \sqrt[3]{-2592 \mu ^4+5670 \mu ^3-5022 \mu ^2+A_9+810 \mu +162}.$\\

\end{document}